\documentclass[12pt]{amsart}
\usepackage{amssymb}
\usepackage{verbatim}

\theoremstyle{plain}
\newtheorem{theorem}{Theorem}
\newtheorem{corollary}{Corollary}
\newtheorem{proposition}{Proposition}

\theoremstyle{definition}
\newtheorem{definition}{Definition}

\theoremstyle{remark}
\newtheorem{remark}{Remark}
\newtheorem{example}{Example}

\newcommand{\cstar}{\ensuremath{\text{C}^{*}}\nobreakdash-\hspace{0 pt}}
\renewcommand{\star}{\ensuremath{{}^{*}}\nobreakdash-\hspace{0 pt}}

\newcommand{\ds}{\displaystyle}

\newcommand{\Alg}{\operatorname{Alg}}
\newcommand{\supp}{\operatorname{supp}}

\newcommand{\com}{\operatorname{com}}
\newcommand{\rad}{\operatorname{rad}}

\newcommand{\ol}[1]{\overline{#1}}
\newcommand{\Uab}{U_{\alpha,\beta}}
\newcommand{\barUab}{U_{\ol{\alpha},\ol{\beta}}}
\newcommand{\Ucd}{U_{\gamma,\delta}}
\newcommand{\Sab}{S^{\phantom{*}}_{\alpha}S^*_{\beta}}
\newcommand{\Scd}{S^{\phantom{*}}_{\gamma}S^*_{\delta}}
\newcommand{\Xab}{\chi_{\alpha,\beta}}
\newcommand{\barXab}{\chi_{\ol{\alpha},\ol{\beta}}}
\newcommand{\barXaa}{\chi_{\ol{\alpha},\ol{\alpha}}}
\newcommand{\barXbb}{\chi_{\ol{\beta},\ol{\beta}}}
\newcommand{\Xcd}{\chi_{\gamma,\delta}}
\newcommand{\Uaa}{U_{\alpha,\alpha}}

\def\bbC{\mathbb C}
\def\bbD{\mathbb D}

\def\bbN{\mathbb N}

\def\bbR{\mathbb R}
\def\bbT{\mathbb T}
\def\bbZ{\mathbb Z}

     \newcommand{\sO}{\mathcal O}

     \newcommand{\sT}{\mathcal T}
     
     \newcommand{\sV}{\mathcal V}

\begin{document}
\title[Cuntz Subalgebras]{Subalgebras of the Cuntz \cstar algebra}
\author{Alan Hopenwasser}
\address{Department of Mathematics\\
        University of Alabama\\
        Tuscaloosa, AL 35487}
\email{ahopenwa@euler.math.ua.edu}

 \author{Justin R. Peters}
 \address{Department of Mathematics\\
          Iowa State University\\
        Ames, IA 50011}
 \email{peters@iastate.edu}

 \date{March 26, 2003}

\thanks{The authors would like to thank Steve Power for helpful
  comments on the subject matter of this paper and for making a major
  contribution to Theorem~\ref{T:spec}.}
 \begin{abstract}
In this paper we exploit the fact that a Cuntz
C$^*$-algebra  is a groupoid C$^*$-algebra 
to facilitate the study of
non-self-adjoint subalgebras of $O_n$.  The Cuntz groupoid is not
principal and the spectral theorem for bimodules does not apply in
full generality.  We characterize the bimodules (over a natural masa)
which are determined by their spectra in the Cuntz groupoid; these are
exactly the ones which are invariant under the guage automorphisms and
exactly the ones which are generated by the Cuntz partial isometries
which they contain.

We investigate analytic subalgebras of $O_n,\ n$ finite, 
by studying cocycles on
the Cuntz groupoid.  In contrast to AF groupoids, there are no
cocycles which are integer valued or bounded and vanish precisely on
the natural diagonal.  $O_n$ contains a canonical UHF subalgebra;
each strongly maximal triangular subalgebra 
of the UHF subalgebra has an
extension to a strongly maximal triangular subalgebra of $O_n$ 
 and
each trivially
analytic subalgebra of the UHF subalgebra has a
proper analytic extension. 

We also study the Volterra subalgebra of $O_n$.  We identify the
spectrum of the Volterra subalgebra and use this to prove a theorem of
Power that the radical is equal to the closed commutator ideal.  We
also show that the Volterra subalgebra is maximal triangular but not
strongly maximal triangular.
 \end{abstract}
\maketitle

\section{Introduction} \label{S:intro}
Non-self-adjoint subalgebras of AF \cstar algebras now have an
extensive theory (see~\cite{MR94g:46001} for an 
introduction); subalgebras of
other classes of \cstar algebras are less well known.
In 1985 Power~\cite{MR86d:47057} investigated the Volterra
subalgebra of a Cuntz \cstar algebra, a topic we shall revisit in this
paper using groupoid techniques.  
More recently, Popescu made use of both the norm closed and the weak
operator topology closed algebras generated by $n$ isometries with
orthogonal range in his multi-variable version of the von Neumann
inequality for a contraction~\cite{MR92k:47073}. 
 And in~\cite{MR2000k:47005}, Davidson and Pitts 
introduced the systematic
study of free semigroup algebras (the WOT closed algebra generated by
$n$ isometries with pairwise orthogonal ranges).  
See~\cite{MR2002a:47107} for one continuation 
of this theme.  When the
range projections of the isometries add to the identity, the free
semigroup algebra is a weakly closed subalgebra of the weak
closure of some representation of $O_n$; these algebras play an
important role in the classification of certain representations of the
Cuntz algebra.

By contrast, the subalgebras which we study are only norm closed.
 More importantly, these subalgebras  are bimodules over
a natural canonical masa in $O_n$. 
(The non-self-adjoint norm closed algebra generated by the generators
of $O_n$ is not such a bimodule.)
It is the bimodule property which
 makes our study amenable to groupoid techniques, so both our
methods and our results have little in common 
with~\cite{MR2002a:47107,MR2000k:47005}. Still, this paper 
and~\cite{MR2002a:47107,MR2000k:47005}
are part of a growing body of work on 
non-self-adjoint algebras related to the Cuntz and
 Cuntz-Krieger algebras.

The Cuntz \cstar algebra $O_n$ is a groupoid 
\cstar algebra~\cite{MR82h:46075,MR2001a:22003}; the Cuntz groupoid
has some similarities to (and one major difference from) the AF
groupoids which determine AF \cstar algebras.  In the spirit of much
of the work on subalgebras of AF \cstar algebras, this paper uses
groupoid techniques to investigate
subalgebras of $O_n$.  The two principal topics
studied are the Volterra subalgebra of $O_n$ (we recover many of
Power's results with simpler proofs as well as a few new facts) and
analytic subalgebras of $O_n$.  Analytic subalgebras of groupoid 
\cstar algebras are most conveniently defined in terms of cocycles on
the groupoid (at the cost of obscuring connections with analyticity);
accordingly, our study of analytic subalgebras of $O_n$ focuses
primarily on cocycles on the Cuntz groupoid. In this paper,
 $n$ will always be finite.

Further on in the introduction we will give a brief review of those
aspects of the Cuntz \cstar algebra which we need as well as a sketch
of the relevant terminology concerning groupoids.  In 
section~\ref{S:groupoid} we will give a more detailed review of the
Cuntz groupoid and its connection with $O_n$.  In
section~\ref{S:spec_th} we describe how to obtain a number of
interesting subalgebras of $O_n$ through the use of the Cuntz
groupoid.  It is here that the most striking difference from AF
algebras appears.  AF groupoids are $r$-discrete principal groupoids
and so the spectral theorem for bimodules~\cite{MR90m:46098} is valid;
this establishes a one-to-one correspondence between bimodules over a
canonical masa and open subsets of the groupoid.  As we show by
example, this theorem is not valid in full generality for the Cuntz
groupoid.  In Theorem~\ref{T:spec}
  we characterize which bimodules over a
natural masa in $O_n$ do arise from open subsets of the Cuntz
groupoid.  

In section~\ref{S:analytic} we discuss cocycles on the Cuntz groupoid
and their associated analytic algebras.  The main tool is
Theorem~\ref{T:cocyc_fun}, which establishes a one-to-one
correspondence between cocycles and continuous real valued functions
defined on the space of units of the groupoid.  
Theorem~\ref{T:nococycles} illustrates further differences from the AF
context.  The most important cocycles used to define analytic
subalgebras of AF \cstar algebras are either integer valued or
bounded.  We show, in contrast, that the Cuntz groupoid supports
neither integer valued nor bounded cocycles which vanish precisely on
the space of units.

The Cuntz \cstar algebra contains an $n^{\infty}$-UHF algebra in a
natural way.  We show in section~\ref{S:analytic} how any 
strongly maximal triangular subalgebra of the 
$n^{\infty}$-UHF subalgebra (with natural diagonal) is contained in a
strongly maximal triangular subalgebra of $O_n$ (with the same
diagonal).   Furthermore,
the refinement TAF subalgebra of the canonical UHF subalgebra
(or any trivially analytic
subalgebra, for that matter) has an extension to an analytic
subalgebra of $O_n$.

Section~\ref{S:Volterra} is devoted to the 
Volterra subalgebra of $O_n$.
We define this algebra through its spectrum in the Cuntz groupoid and
then show that it is equal to the intersection of a natural
representation of $O_n$ acting on $L^2[0,1]$ and the usual Volterra
nest subalgebra of $L^2[0,1]$.
We identify the spectrum of the radical of the Volterra nest
subalgebra of $O_n$ and use this to prove Power's theorem that the
radical is equal to the closed commutator ideal. 
The Volterra subalgebra provides
an example of a maximal triangular subalgebra of $O_n$ which is not
strongly maximal.

We conclude the introduction with some salient facts about $O_n$ and
some terminology concerning groupoids.  Cuntz proved 
in~\cite{MR57:7189} that, up to isomorphism, there is only one 
\cstar algebra generated by $n$ isometries whose range projections are
pairwise orthogonal and add to the identity operator.  If
$S_1, S_2, \dots, S_n$ are generators of $O_n$, then 
$S^*_iS^{\vphantom{*}}_j$ is either 0 (if $i \ne j$) or $I$
(if $i=j$).  Using this, it is easy to see that any word in the
$S_i$ and the $S^*_i$ can be written in the form
$S_{\alpha_1}S_{\alpha_2}\cdots S_{\alpha_k}
S^*_{\beta_j}\cdots S^*_{\beta_1}$ for some finite strings
$\alpha = (\alpha_1, \dots, \alpha_k)$ and
$\beta = (\beta_1, \dots, \beta_j)$ whose coordinates are taken from
$\{1, 2, \dots, n\}$.

If $\alpha$ and $\beta$ are finite
strings, let $\Sab$ denote the partial isometry
$S_{\alpha_1}S_{\alpha_2}\cdots S_{\alpha_k}
S^*_{\beta_j}\cdots S^*_{\beta_1}$.  Either string,
$\alpha$ or $\beta$ could be empty.  If both are empty,
$\Sab = I$.
The set of all $\Sab$ forms an inverse semigroup
(the \emph{Cuntz inverse semigroup}
generated by $S_1, \dots S_n$).  Each element of the Cuntz inverse
semigroup will be referred to as a \emph{Cuntz partial isometry}.

A particularly useful representation of the Cuntz \cstar algebra
acts on $L^2[0,1]$.  If we partition
$[0,1]$ into $n$ disjoint subintervals (ignoring endpoints), then the
order preserving
affine map from $[0,1]$ onto the $i^{\text{th}}$ subinterval
 induces an isometry $S_i$ on $L^2[0,1]$.  The $S_i$ generate 
a copy of the Cuntz algebra.  

Each partial isometry $\Sab$ is induced by an
order preserving  affine
map from some subinterval of $[0,1]$ onto some other subinterval; we
often
let these affine maps ``stand in'' for the partial isometries when
discussing them.  So, for example, $S_i$ is the partial isometry from
$[0,1]$ onto the $i^{\text{th}}$ subinterval of $[0,1]$ -- using the
natural order of the indices to match the natural order of the
subintervals.  In this fashion, $S_{(i,j)}=S_iS_j$
is the partial isometry from
$[0,1]$ onto the $j^{\text{th}}$ subinterval of the $i^{\text{th}}$
subinterval of $[0,1]$. (To be specific about endpoints,
the $i^{\text{th}}$ subinterval of $[0,1]$ is
$[\frac {i-1}n, \frac in]$ and the $j^{\text{th}}$ subinterval of
this is 
$[\frac {i-1}n + \frac {j-1}{n^2}, \frac {i-1}n + \frac j{n^2}]$.)
 This latter ``range'' interval has length
$1/n^2$.  The adjoint, $S^*_{(i,j)} = S_j^*S_i^*$ is, of course, the
partial isometry from the $j^{\text{th}}$ subinterval of
 the $i^{\text{th}}$ subinterval of $[0,1]$ onto  $[0,1]$.
Thus, the range interval for $S_{\alpha}$ is determined by looking at
interval $\alpha_1$ in the partition of $[0,1]$ into $n$ parts, 
then at interval $\alpha_2$ in the partition of
interval $\alpha_1$ into $n$ parts, 
then at interval $\alpha_3$ in the partition of 
the preceding interval into $n$ parts, etc.
 The domain of $S^*_{\alpha}$ is, of
course, the same as the range of $S_{\alpha}$.  Based on this model,
(and using the notation $\ell(\alpha)$ to denote the
\emph{length} of $\alpha$, the number of coordinates in the string)
if $\ell(\alpha) \geq \ell(\beta)$ we refer to $\Sab$
as being \emph{contractive}
(the range interval is shorter than the domain
interval); if $\ell(\alpha) \leq \ell(\beta)$
we say that $\Sab$ is \emph{expansive}.

The Volterra subalgebra of $L^2[0,1]$ is defined as the algebra of all
operators which leave invariant each projection onto
$L^2[0,x]$ viewed as a subspace of $L^2[0,1]$.  Let $p_x$ denote this
projection and $\sV$ the nest of all such projections.  Most of the
$p_x$'s do not lie in $O_n$; in fact, $p_x \in O_n$ if, and only if,
$x$ is an $n$-adic rational in $[0,1]$.  However,
$\{p_x \mid x \text{ is $n$-adic}\}$ is strongly dense in $\sV$, so
either nest can be used to determine the Volterra nest algebra.

Finally, we briefly review some groupoid terminology.
A groupoid is a set $G$ with a
partially defined multiplication and an inversion.  If $a$ and $b$ can
be multiplied, then $(a,b)$ is called a \emph{composable pair};
$G^2$ denotes the set of all composable pairs.  The multiplication
satisfies an associative law and inversion satisfies
$(a^{-1})^{-1}=a$.  Elements of the form $a^{-1}a$ and $aa^{-1}$ are
called \emph{units}. Units act as left and right identities when
multiplied by elements with which they are composable.
Groupoids have \emph{range} and
\emph{domain} maps defined by $r(a) = aa^{-1}$ and
$d(a) = a^{-1}a$.  A subset $E \subseteq G$ is said to be a 
$G$-\emph{set} if $r$ and $d$ are both one-to-one on $E$.

The Cuntz groupoid $G_n$ (described in detail in 
section~\ref{S:groupoid}) is a locally compact space whose topology is
generated by a collection of compact open $G$-sets.  The unit space is
open, so $G_n$ is $r$-discrete (by definition); but $G_n$ is not a
principal groupoid (i.e., it is not an equivalence relation on the
unit space).

The groupoid \cstar algebra of $G_n$ is built from
$C_c(G_n)$, the set of continuous functions on $G_n$ with compact
support.  A convolution type multiplication and an inversion turn 
$C_c(G_n)$ into a \star algebra; the definitions of these operations
specific to this context are given in section~\ref{S:groupoid}.  An
auxiliary norm $\|\ \|_I$ is put on $C_c(G_n)$ and a \cstar norm is defined on 
$C_c(G_n)$ by $\|f\| = \sup_{\pi} \|\pi(f)\|$, where $\pi$ runs over
all \star representations of $C_c(G_n)$ which are decreasing with
respect to $\|\ \|_I$.  The groupoid \cstar algebra $C^*(G_n)$ is the
completion of $C_c(G_n)$ with respect to the \cstar norm.

Extensive treatments of groupoids and their associated \cstar algebras
can be found in~\cite{MR82h:46075} and~\cite{MR2001a:22003}.  We refer
the reader to either of those monographs for details.  Renault 
proves~\cite[Proposition 4.1]{MR82h:46075} that for $r$-discrete
groupoids, the $\|\ \|_{\infty}$-norm 
on $C_c(G)$ is dominated by the 
\cstar norm.  As a consequence, any element of $C^*(G)$ can be
identified with a continuous function on $G_n$ which vanishes at
infinity.  Caveat: not every $C_0$ function on $G$ is associated with
an element of the groupoid \cstar algebra.  We shall make extensive
use of this identification throughout this paper; in particular,
elements of $C^*(G_n)$ will routinely be viewed as continuous functions
on $G_n$ (vanishing at infinity).  
Both~\cite{MR82h:46075} and~\cite{MR2001a:22003} identify $C^*(G_n)$
as $O_n$; all the same, for the convenience of the reader we indicate
the connection between $C^*(G_n)$ and $O_n$ in
section~\ref{S:groupoid}. 

\section{Review of the Cuntz groupoid} \label{S:groupoid}

The following description of the groupoid $G_n$ for the Cuntz algebra,
$O_n$ is taken from Paterson's book~\cite{MR2001a:22003}.
(An earlier description of
this groupoid is in Renault's book~\cite{MR82h:46075}.)
Much of what we outline in this section appears in more general
form in~\cite{MR98g:46083}.  The material here is sufficient for
subsequent sections of this paper and should prove convenient for
readers unfamiliar with the Cuntz groupoid.

Let $\ds X = \prod_1^{\infty}\{1,2,\dots,n\}$; equip $X$ with the
product topology.

For a finite string $\alpha$ of digits in $\{1,2,\dots,n\}$, let
$\ell(\alpha)$ denote the length of the string.  If $\alpha$ is a
finite sequence and $\gamma$ is an element of $X$, then
$\alpha\gamma$ denotes the element of $X$ obtained by `prepending'
$\alpha$ to $\gamma$.  As a set, the Cuntz groupoid is
\[
G_n = \{(\alpha\gamma, \ell(\alpha)- \ell(\beta), \beta\gamma)
\mid \alpha \text{ and } \beta \text{ are finite strings and }
\gamma \in X\}.
\]
Two elements $(x,k,y)$ and $(w,j,z)$ of $G_n$ are composable if, and
only if, $y=w$.
Multiplication and inversion in the groupoid are given by
the formulas:
\[
(x,k,y)(y,j,z) = (x,k+j,z) \quad \text{and} \quad
(x,k,y)^{-1} = (y,-k,x).
\]

Multiplication in $C^*(G_n)$ is given by a
convolution formula:
\[
f \ast g (x,k,y) = \sum f(x,k_1,z)g(z,k_2,y)
\]
The sum is taken over all $k_1, k_2 \in \bbZ$ and $z \in X$ such that
$(x,k_1,z)$ and $(z,k_2,y)$ are in $G_n$ and $k_1+k_2=k$.  Initially,
this formula is used in $C_c(G_n)$; in this context all but finitely
many terms in the sum are zero.  But the same formula is also valid
for all those functions in $C_0(G_n)$ which correspond to elements of
$C^*(G_n)$.   The adjoint is given by the formula
\[
f^*(x,k,y) = \ol{f(y,-k,x)}
\]

The Cuntz groupoid is a locally compact groupoid.
 The sets
\[
U_{\alpha,\beta} = \{(\alpha\gamma, \ell(\alpha) - \ell(\beta),
\beta\gamma \mid \gamma \in X \}.
\]
form a basis for the topology.

A neighborhood basis at a 
point $(x,k,y) \in G_n$ can be obtained in
the following way.  Write $x = (x_1, x_2, x_3, \dots)$ and
$y=(y_1,y_2,y_3,\dots)$.  There is a positive integer $p$ such that
$y_j = x_{j+k}$ for all $j \geq p$.  We may as well assume that $p$ is
the smallest integer which satisfies this property (although this is
not essential).  When $k<0$, we also assume that $-k<p$, so that
$j+k \geq 1$ when $j \geq p$.
  For each $j \geq p$, let
\begin{align*}
\alpha^j &= (x_1, \dots, x_{j+k}), \\
\beta^j &= (y_1, \dots, y_j).
\end{align*}
Then $\ell(\alpha^j) - \ell(\beta^j) = j+k-j=k$.  If
\[
\gamma^j = (x_{j+k+1}, x_{j+k+2}, \dots) =
(y_{j+1}, y_{j+2}, \dots ),
\]
then $(x,k,y) = (\alpha^j\gamma^j, \ell(\alpha^j) - \ell(\beta^j),
\beta^j\gamma^j)$.  So, $(x,k,y) \in U_{\alpha^j,\beta^j}$, for
all $j \geq p$.

On the other hand, let
$(\ol{x}, \ol{k}, \ol{y}) \in G_n$ and suppose that
$\ol{x} \neq x$.  Then, for some $i$,
$\ol{x}_i \neq x_i$. Consequently, for all $j > i$,
$(\ol{x}, \ol{k}, \ol{y}) \notin
U_{\alpha^j,\beta^j}$.  Similarly, if
$\ol{y} \neq y$ or, even more easily, if
$\ol{k} \neq k$, then
$(\ol{x}, \ol{k}, \ol{y}) \notin
U_{\alpha^j,\beta^j}$, for large $j$.  The conclusion is that
\[
\{(x,k,y)\} = \cap_{j=p}^{\infty} U_{\alpha^j,\beta^j};
\]
 so $\{ U_{\alpha^j,\beta^j} \mid j \geq p \}$ forms a neighborhood
basis at $(x,k,y)$.  

If we restrict the topology to a basic open set $\Uab$, 
  the association
\[
\gamma \longleftrightarrow
(\alpha\gamma, \ell(\alpha) - \ell(\beta), \beta\gamma)
\]
is a homeomorphism between $X$ and $U_{\alpha,\beta}$.  In particular,
each $U_{\alpha,\beta}$ is compact in the topology on $G_n$.

The groupoid \cstar algebra of the Cuntz groupoid $G_n$ is the Cuntz 
\cstar algebra $O_n$~\cite{MR82h:46075,MR2001a:22003}.  The rest of
this section is intended to illuminate the connection between the
groupoid and the \cstar algebra.

Fix
\begin{align*}
\alpha &= (\alpha_1, \dots, \alpha_a), \\
\beta &= (\beta_1, \dots, \beta_b), \\
\gamma &= (\gamma_1, \dots, \gamma_c), \\
\delta &= (\delta_1, \dots, \delta_d).
\end{align*}
Suppose $b<c$ and let
$\ol{\gamma} = (\gamma_{b+1},\dots,\gamma_c)$.
Since
\[
S_i^*S_j =
\begin{cases}
0, &\text{if $i \neq j$,} \\
I, &\text{if $i=j$,}
\end{cases}
\]
we have
\[
\Sab \Scd =
\begin{cases}
S_{\alpha}S_{\ol{\gamma}}S^*_{\delta},
&\text{if $(\beta_1, \dots, \beta_b) = (\gamma_1, \dots, \gamma_b)$},
  \\
0, &\text{if $(\beta_1, \dots, \beta_b) \neq (\gamma_1, \dots, \gamma_b)$}.
\end{cases}
\]

In the following, $\Xab$ will denote the characteristic
function of $\Uab$.
  This is a continuous function
on $G_n$ with compact support.
The convolution, $\Xab \ast \Xcd$ is given
by the formula
\[
\Xab \ast \Xcd (x,k,y) =
\sum \Xab (x,i,u) \Xcd (u,j,y),
\]
where $i+j=k$ and $u$, $i$, and $j$ run through the countably many
choices yielding composable elements of $G_n$.

Suppose  $u$ is such that $\Xab (x,i,u) \Xcd (u,j,y) \neq 0$.  Then
\begin{align*}
(u_1, \dots, u_b) &= (\beta_1, \dots, \beta_b) \text{ and} \\
(u_1, \dots, u_c) &= (\gamma_1, \dots, \gamma_c).
\end{align*}
With the assumption $b < c$ in force, this yields,
$(\beta_1, \dots, \beta_b) = (\gamma_1, \dots, \gamma_b)$.

Thus, if
$(\beta_1, \dots, \beta_b) \neq (\gamma_1, \dots, \gamma_b)$, then
every term in the convolution sum is zero and we have the
correspondence
\[
0=\Xab \ast \Xcd \longleftrightarrow \Sab
\Scd=0.
\]

So  assume that
$(\beta_1, \dots, \beta_b) = (\gamma_1, \dots, \gamma_b)$.  Recall
$\ol{\gamma} = (\gamma_{b+1}, \dots, \gamma_c)$ and
$\ell(\ol{\gamma}) = c-b$.
Suppose $i$, $j$ and $u$ are such that
$\Xab (x,i,u) \Xcd (u,j,y) \neq 0$.  This forces
\begin{align*}
i &= \ell(\alpha) - \ell(\beta) = a-b, \text{ and} \\
j &= \ell(\gamma) - \ell(\delta) = c-d.
\end{align*}
and hence
\[
i+j = a+c-b-d = \ell(\alpha) + \ell(\ol{\gamma})
- \ell(\delta) = \ell(\alpha \ol{\gamma}) -\ell(\delta).
\]
Also, as before,
\begin{align*}
(u_1, \dots, u_b) &= (\beta_1, \dots, \beta_b), \text{ and}\\
(u_1, \dots, u_c) &= (\gamma_1, \dots, \gamma_c).
\end{align*}
In addition,
\begin{align*}
(x_{a+1}, x_{a+2}, \dots) &= (u_{b+1}, u_{b+2}, \dots),
\text{ and} \\
(u_{c+1}, u_{c+2}, \dots) &= (y_{d+1}, y_{d+2}, \dots).
\end{align*}
Hence,
\[
(x_{a+1}, \dots, x_{a+c-b}) = (u_{b+1}, \dots, u_c)
=(\gamma_{b+1}, \dots, \gamma_c) = \ol{\gamma}.
\]
Letting $\eta = (x_{a+c-b+1}, \dots)$, we see that
$x = \alpha \ol{\gamma} \eta$.
Now
\[
(y_1, \dots, y_d) = (\delta_1, \dots, \delta_d) = \delta
\]
and
\[
\eta = (x_{a+c-b+1}, \dots) = (u_{b+c-b+1}, \dots)
=(u_{c+1}, \dots) = (y_{d+1}, \dots),
\]
so $y = \delta \eta$.

Since $u$ is completely determined by these conditions, there is only
one term in
\[
\sum \Xab (x,i,u) \Xcd (u,j,y)
\]
which is non-zero if the sum is non-zero.  This happens only if
$x = \alpha \ol{\gamma} \eta$ and $y = \delta \eta$.

Conversely, if $x=\alpha \ol{\gamma} \eta$ and
$y = \delta \eta$, for some $\eta$ in $X$, then choose $u$ to be
$\beta \ol{\gamma} \eta$; for this choice of $u$
\[
\Xab (x, a-b, u) \Xcd (u, c-d, y) = 1.
\]
(No other choice of $u$ yields a non-zero value for this
product.)  Thus,
\[
\Xab \ast \Xcd (x,k,y) = 1 \Longleftrightarrow
\begin{cases}
k &= \ell(\alpha \ol{\gamma}) - \ell(\delta), \\
x &= \alpha \ol{\gamma} \eta, \\
y &= \delta \eta, \qquad \text{for some $\eta \in X$.}
\end{cases}
\]
Thus,
\[
\Xab \ast \Xcd (x,k,y) = 1 \Longleftrightarrow
\chi_{\alpha \ol{\gamma}, \delta}(x,k,y) = 1
\]
and
\[
\Xab \ast \Xcd (x,k,y) = 0 \Longleftrightarrow
\chi_{\alpha \ol{\gamma}, \delta} (x,k,y)= 0.
\]
The only possible values for $\Xab \ast \Xcd$ are $0$ and $1$, so
$\Xab \ast \Xcd = \chi_{\alpha \ol{\gamma}, \delta}$.  Thus we have
  the correspondence
\[
\Xab \ast \Xcd \longleftrightarrow
S_{\alpha}S_{\ol{\gamma}}S^*_{\delta} =
S_{\alpha} S^*_{\beta} \Scd.
\]

In the case in which $b=c$, omit the $\ol{\gamma}$.  If $b>c$, a
similar argument yields the correspondence between
$\Xab \ast \Xcd$ and
$S_{\alpha} S^*_{\beta} S_{\gamma} S^*_{\delta}$.
It is now clear that 
$\chi_{1,\emptyset}, \dots, \chi_{n,\emptyset}$ generate all
$\Xab$ in $C^*(G_n)$ and that
$\chi_{1,\emptyset}, \dots, \chi_{n,\emptyset}$ are isometries whose
range projections sum to $I$.  This illuminates the identification of
$C^*(G_n)$ as $O_n$.

\section{The spectral theorem for bimodules} \label{S:spec_th}

In this section we show how to associate subalgebras of $O_n$ to
certain subsets of the Cuntz groupoid.
The spectral theorem for bimodules, proven in~\cite{MR90m:46098}
for principal groupoids, is
extremely useful in the study of non-self adjoint subalgebras of
groupoid \cstar algebras;  however, this theorem is not valid
in full generality for the Cuntz groupoid.  
In Theorem~\ref{T:spec} we characterize those bimodules over a natural
masa for which the spectral theorem is valid. In what follows we assume all bimodules
are norm closed.

Keep in mind in the sequel that we view elements of
$C^*(G_n) = O_n$ as continuous functions on $G_n$ which vanish at
infinity.  Also, let $D_0$ denote the space of units in $G_n$;
i.e., $D_0 = \{(x,0,x) \mid x \in X \}$.  This is an open subset of
$G_n$; in fact, it is $\Uab$, where both $\alpha$ and $\beta$ are
empty strings.

If $P$ is an open subset of $G_n$, let
\[
A(P) = \{f\in \text{C}^{*}(G_n) \mid \supp f \subseteq P\}.
\]
If $f \in A(P)$, then $f(x,k,y)=0$ for all
$(x,k,y) \in G_n \setminus P$.
 For any open set
$P$ in $G_n$, $A(P)$ is a bimodule over $A(D_0)$.  This is easy to
see:   when $f \in A(D_0)$, $f$ is
supported on $D_0$, so
\begin{align*}
f \ast g (x,k,y) &= f(x,0,x)g(x,k,y) \\
\intertext{and}
g \ast f (x,k,y) &= g(x,k,y)f(y,0,y).
\end{align*}
If $g$ is supported in $P$ and $f \in A(D_0)$,
 then $f \ast g$ and $g \ast f$ are also supported in $P$.

If, in addition, $P$ satisfies the property
\[
(x,k,y) \in P \text{ and } (y,j,z) \in P \Longrightarrow
(x,k+j,z) \in P,
\]
then $A(P)$ is a subalgebra of $O_n$. 
  Also, if $D_0 \subseteq P$, then
$A(D_0) \subseteq A(P)$.

One choice for $P$ is
$D_0$.  We shall reprove in
Proposition~\ref{cmasa} the known fact
that $A(D_0)$ is a masa in $O_n$. 
Another choice for $P$ is
\[
P_{UHF} = \{(x,k,y) \in G_n \mid k=0 \}.
\]
Then $A(P_{UHF})$ is a subalgebra of $O_n$ which contains
$A(D_0)$; the first part of Theorem~\ref{T:spec} shows that
$A(P_{UHF})$ is generated by the Cuntz partial isometries which it
contains.  But these are just the $\Sab$ with 
$\ell(\alpha) = \ell(\beta)$; it is well known that the subalgebra of
$O_n$ generated by these Cuntz partial isometries is an
$n^{\infty}$-UHF algebra.  We shall refer to this subalgebra as the
\emph{canonical UHF subalgebra} of $O_n$.  Note also that there is an
obvious isomorphism between $P_{UHF}$ and the usual groupoid for the
$n^{\infty}$-UHF algebra.

If $B \subseteq O_n$ is an $A(D_0)$ bimodule, define
\[
\sigma(B) = \{(x,k,y) \mid \text{there is } f \in B
\text{ with } f(x,k,y) \ne 0 \}.
\]
In other words,
$G_n \setminus \sigma(B) = \{(x,k,y) \mid f(x,k,y)=0 \text{ for all }
f \in B \}$. Clearly, $\sigma(B)$ is an open subset of $G_n$.

\begin{proposition} \label{P:spec}
If $P$ is an open subset of $G_n$, then 
$\sigma(A(P)) = P$. 
\end{proposition}

\begin{proof}
Suppose that  $P$ is an open subset of $G_n$.
If $(x,k,y) \notin P$, then
$f(x,k,y) = 0$ for all $f \in A(P)$, so
$(x,k,y) \in G_n \setminus \sigma(A(P))$.  Thus
$\sigma(A(P)) \subseteq P$.  On the other hand, suppose that
$(x,k,y) \in P$.  Since $P$ is open, there are $\alpha$ and
$\beta$ such that $\ell(\alpha) - \ell(\beta) =k$ and
$(x,k,y) \in \Uab \subseteq P$. 
Now $\Xab \in C^*(G_n)$
(and corresponds to the partial isometry $\Sab$).  Since
$\supp \Xab \subseteq P$, $\Xab \in A(P)$.  But
$\Xab (x,k,y) =1$, so $(x,k,y) \in \sigma(A(P))$.  This shows that
$P \subseteq \sigma(A(P))$.
\end{proof}

\begin{definition}
A bimodule over $A(D_0)$ is said to be \emph{reflexive} if
$B = A(\sigma(B))$.
\end{definition}

Note that the reflexive bimodules over $A(D_0)$ are exactly the ones
of the form $A(P)$ for some open set $P \subseteq G_n$.
If the spectral theorem for bimodules over $A(D_0)$ were valid in full
generality, then every bimodule over $A(D_0)$ would be reflexive.
The following example gives a bimodule which is not reflexive.
The authors thank Steve Power for drawing this
example to their attention.

\begin{example}
 Let $P = D_0 \cup U_{(1),\emptyset}$.  
Note that $D_0$ is the support set for $I$ and 
$ U_{(1),\emptyset}$ is the support set for $S_1$.
Let $\ol{1}$ denote the element $(1,1,1,\dots)$ in $X$.
Both $(\ol{1},0,\ol{1})$ and $(\ol{1},1,\ol{1})$ are elements
of $P$.  Let
\[
B = \{f \in O_n \mid \supp f \subseteq P \text{ and }
f(\ol{1},0,\ol{1}) = f(\ol{1},1,\ol{1})\}.
\]
If $f \in B$ and $g \in A(D_0)$, then
\begin{align*}
g \ast f (\ol{1},0,\ol{1}) 
&= g(\ol{1},0,\ol{1}) f(\ol{1},0,\ol{1})\\
&= g(\ol{1},0,\ol{1}) f(\ol{1},1,\ol{1})\\
&= g \ast f (\ol{1},1,\ol{1}).
\end{align*}
Similarly $f \ast g (\ol{1},0,\ol{1}) =
f \ast g (\ol{1},1,\ol{1})$; thus $B$ is a bimodule.
From the definition of $B$ it is clear that 
$\sigma(B) \subseteq P$.  But the characteristic function of $P$ is in
$B$, so $P \subseteq \sigma(B)$.  Thus $\sigma(B) = P$ and it is
trivial that $B \ne A(P)$.

Note that $B$ is the bimodule over $A(D_0)$ generated by
$I + S_1$.
\end{example}

The next theorem characterizes the reflexive bimodules over 
$A(D_0)$.  One of the conditions equivalent to reflexivity is
invariance under the guage automorphisms.  For each complex number 
$\lambda$ of absolute value one, a guage automorphism 
of $O_n$ is determined by
its action on the generators:
$\eta_{\lambda}(S_i) = \lambda S_i$.  
The authors thank Steve Power for pointing out 
condition~(\ref{guage_inv}) and providing a proof
that~(\ref{guage_inv}) implies~(\ref{gen_Cuntz}).

\begin{theorem}[Spectral Theorem for Bimodules] \label{T:spec}
Let $B$ be a bimodule over $A(D_0)$.  Then the following are
equivalent:
\begin{enumerate}
\item  \label{reflexive} $B$ is reflexive.
\item \label{gen_Cuntz}  $B$ is generated by the Cuntz
partial isometries which it contains. 
\item  \label{guage_inv} $B$ is invariant under 
all the guage automorphisms.
\end{enumerate}
\end{theorem}

\begin{proof}
We first show~(\ref{reflexive}) implies~(\ref{gen_Cuntz}); i.e., 
 any bimodule of the form $A(P)$ with $P$ an open
subset of $G_n$ is generated by the Cuntz partial isometries which it
contains.  Let $B$ be the bimodule generated by the Cuntz partial
isometries in $A(P)$.  Suppose that $\Xab$ is a Cuntz partial isometry
in $A(P)$, so that $\Uab \subseteq P$.  Let $f$ be any continuous
function supported on $\Uab$.  Define a function $g$ on $\Uaa$ by 
$g(\alpha\gamma,0,\alpha\gamma) =
f(\alpha\gamma, \ell(\alpha)-\ell(\beta), \beta\gamma)$.  
Observe that $f = g \ast \Xab$.  Since $g \in A(D_0)$, we have
$f \in B$.  Thus, $B$ contains any continuous function supported on a
compact open set of the form $\Uab \subseteq P$.  Since any compact
open subset of $P$ can be written as a finite union of sets of the
form $\Uab$, $B$ contains any continuous function supported on a
compact open subset of $P$. 
 Any compact set is contained in a compact open set,
so all continuous functions with compact support in $P$
 are in $B$.  But
these are dense (in the \cstar norm) in $A(P)$, so $A(P) = B$.

To prove that~(\ref{gen_Cuntz}) implies~(\ref{reflexive}),
 suppose that $B$ is an $A(D_0)$ bimodule and that $B$ is generated
by the Cuntz partial isometries which it contains.  Let 
$P = \bigcup \Uab$, where the union is taken over all $\alpha,\beta$
such that $\Xab$ is in $B$.  It is obvious that
$P \subseteq \sigma(B)$.  To see the reverse containment, suppose that
$(x,k,y) \notin P$.  Then $\Xab (x,k,y) = 0$ for all Cuntz partial
isometries $\Xab$ in $B$.  If $f$ and $g$ are any elements of
$A(D_0)$, then $f \ast \Xab \ast g$ also vanishes at $(x,k,y)$, since 
$f \ast \Xab \ast g (x,k,y) = f(x,0,x) \Xab(x,k,y) g(y,0,y)$.
It follows that any element of the bimodule generated by these 
$\Xab$, i.e., any element of $B$, also vanishes at $(x,k,y)$.  Thus
$(x,k,y) \notin \sigma(B)$.

We now know that $P = \sigma(B)$, so $B \subseteq A(P)$.  It remains
to show that $A(P) =B$.  If $\alpha$ and $\ol{\alpha}$ are
finite strings with 
$\ell(\alpha) < \ell(\ol{\alpha})$ and
$\alpha_j = \ol{\alpha}_j$ for $j = 1, \dots, \ell(\alpha)$, we
shall call $\ol{\alpha}$ an \emph{extension} of $\alpha$.
If $\ol{\alpha}$ is an extension of $\alpha$ and 
$\ol{\beta}$ is an extension of $\beta$ and if
$\ell(\ol{\alpha}) - \ell(\ol{\beta})
= \ell(\alpha) - \ell(\beta)$, then $\barUab \subseteq \Uab$.
Since $\barXab = \barXaa \ast \Xab \ast \barXbb$, if $\Xab$ is in $B$
then so is $\barXab$.

To prove that $A(P) \subseteq B$, it suffices to show that if
$\Xcd$ is a Cuntz partial isometry in $A(P)$, then
$\Xcd$ is in $B$. (We have already seen that $A(P)$ is generated by
the Cuntz partial isometries which it contains.)  Let 
$(x,k,y) \in \Ucd$.  Since $P = \sigma(B)$, there is a Cuntz partial
isometry $\Xab$ in $B$ such that $(x,k,y) \in \Uab$.  For some
extensions $\ol{\alpha}$ of $\alpha$
and $\ol{\beta}$ of $\beta$ with
$\ell(\ol{\alpha}) - \ell(\ol{\beta}) = k$, we have
$(x,k,y) \in \barUab \subseteq \Uab \cap \Ucd$ and, as noted above,
$\barXab \in B$.  

Since $\Ucd$ is compact, we can write it as a finite disjoint union of
sets of the form $\barUab$ with $\barXab \in B$.
  But $\Xcd$ is the sum
of the corresponding $\barXab$; thus $\Xcd \in B$ and we have shown
that $A(P) =B$.

(\ref{gen_Cuntz}) implies (\ref{guage_inv}) is trivial: a guage
automorphism maps a Cuntz partial isometry to a scalar multiple of
itself.   It remains only to prove that (\ref{guage_inv}) implies
(\ref{gen_Cuntz}).  Cuntz showed in~\cite{MR57:7189} that each element
of $O_n$ has a Fourier series with respect to a chosen generator, say
$S_1$:
\[
a \sim \sum_{-\infty}^{-1} (S^*_1)^{|k|}a_k +
\sum_0^{\infty} a_k S^k_1.
\]
The coefficients $a_k$ all lie in the canonical 
UHF subalgebra.  If $\eta_{\lambda}$ is a guage
automorphism, then
\[
\eta_{\lambda}(a) \sim 
\sum_{-\infty}^{-1} (S^*_1)^{|k|} \ol{\lambda}^k a_k +
\sum_0^{\infty} \lambda^k a_k S^k_1.
\]
By Remark 1 of section 1.10 of~\cite{MR57:7189},
there is a Cesaro convergence of generalized polynomials and,
for $k \geq 0$,
\[
a_k S^k_1 = \int_{\bbT} \ol{\lambda}^k 
\eta_{\lambda}(a) d\lambda,
\]
where $d\lambda$ is normalized Lebesgue measure on $\bbT$.
A similar formula holds when $k$ is negative.
Thus, if $a \in B$, then all $a_k S_1^k$ and all
$(S^*_1)^{|k|}a_k$ lie in $B$.
Note that we may, without loss of generality, assume that
each coefficient $a_k$ satisfies
$a_k = a_k S^k_1 (S^*_1)^k$ when $k \geq 0$, with a similar assertion
when $k < 0$.

Let
\[
B_k = \{b \in A(P_{UHF}) \mid bS^k_1 \in B, 
b = bS^k_1 (S^*_1)^k \}
\]
when $k \geq 0$; a similar definition is used when $k<0$.
So $B$ is the closed linear span of all the subspaces
$B_k S^k_1$ and $(S^*_1)^k B_{-k}$.  It is trivial that $B_k$ is a
left bimodule over $A(D_0)$ and easy to see that $B_k$ is a right
bimodule.  (It suffices to prove this for right multiplication 
by projections in $A(D_0)$ and we can limit ourselves to
subprojections of the range projection for $S_1^k$.  But if $p$ is a
subprojection of  $S_1^k (S^*_1)^k$, then there is a subprojection 
$q$ of the domain projection for $S_1^k$  such that
$pS_1^k = S_1^kq$.  For any $b \in B_k$, 
$bpS_1^k = bS_1^kq \in B$ and so $bp \in B_k$.)

Now $B_k$ is a bimodule over $A(D_0)$ contained in the canonical UHF
algebra $A(P_{UHF})$ and the spectral theorem for bimodules is valid in
full generality in this context.  So $B_k$ is generated by the Cuntz
partial isometries which it contains.  It follows that the same is
true for $B$.
\end{proof}

It is shown in~\cite[Remark 2.1.8]{MR82f:46073a} that 
$A(D_0)$ is a masa
in $O_n$; for the convenience of the reader,
we reprove this using groupoid techniques.

\begin{proposition} \label{cmasa}
$A(D_0)$ is a masa in $C^*(G_n) = O_n$.
\end{proposition}

\begin{proof}
Suppose that $g \in O_n$ and
$f \ast g = g \ast f$, for all $f \in A(D_0)$.
We must show that $g \in A(D_0)$.  When $f \in A(D_0)$, $f$ is
supported on $D_0$, so
\begin{align*}
f \ast g (x,k,y) &= f(x,0,x)g(x,k,y) \\
\intertext{and}
g \ast f (x,k,y) &= g(x,k,y)f(y,0,y).
\end{align*}
If $x \neq y$, we can find $f \in A(D_0)$ such that
$f(x,0,x) \neq f(y,0,y)$, hence $g(x,k,y)=0$ when $x \neq y$.
If $k \neq 0$ then $\{(x,k,y)\in G_n  | x=y \}$ has empty interior
and so cannot support a non-zero continuous function.  Thus, if $g$ is
continuous on $G_n$ and $g(x,k,y)=0$ whenever $x \neq y$ and
$k \neq 0$, then $g(x,k,x) = 0$ for all $x$.  

We have shown that if $g$ commutes with $f$ then $g(x,k,y)=0$ whenever
$x \neq y$ or $k \neq 0$.  This means that $g$ is supported on $D_0$
and $A(D_0)$ is a masa in $C^*(G_n)$.
\end{proof}

\section{Analytic subalgebras of $\sO_n$} \label{S:analytic}

\begin{definition}
A \emph{cocycle}  is a continuous 
function $d$ from $G_n$ to $\bbR$
which satisfies
 the  identity
\[ 
d(x, k, y) + d(y, l, z) = d(x, k+l, z) 
\]
for $(x, k, y)$ and  $(y, l, z)$ in $G_n$.

If $P= \{(x,k,y) \mid d(x,k,y) \geq 0 \}$ is open, then we call 
$A(P)$ the \emph{analytic subalgebra} associated with $d$. 
\end{definition}

To avoid trivialities, we always assume that $d$ is not identically
equal to 0.

\begin{remark}
There is a standard procedure for associating a one parameter family
of automorphisms  to a cocycle $d$.  For each 
$t \in \bbR$, an automorphism $\eta_t$ is given by the formula:
\[
(\eta_t f)(x,k,y) = e^{itd(x,k,y)}f(x,k,y). 
\]
Each $\eta_t$ is a \star automorphism of $C_c(G_n)$ onto itself; it is
not hard to show that this automorphism preserves the 
\cstar norm and so extends to an automorphism of $O_n$ with the
formula above.  For each point $(x,k,y) \in G_n$, the map
$f \mapsto f(x,k,y)$ is decreasing with respect to the
$\|\ \|_{\infty}$, and hence also norm decreasing with respect to the 
\cstar norm.  So these maps are continuous linear functionals on
$O_n$.  If $f \in O_n$, we consider all functions of the form
$t \mapsto \rho(\eta_t(f))$, where $\rho$ is a linear functional of
the type above.  Given $f \in O_n$, it is easy to check that
$t \mapsto \rho(\eta_t(f))$ is an $H^{\infty}$-function on
$\bbR$ for all linear functional of this form if, and only if,
$f$ is supported on
$\{(x,k,y) \mid d(x,k,y) \geq 0 \}$.  This indicates why the term
\emph{analytic subalgebra} is appropriate in the definition above.
\end{remark}

\begin{example}
Consider the \emph{dilation cocycle} $d(x,k,y)=k$. 
 Since each Cuntz partial
isometry $\Sab$ corresponds to the function $\Xab$ on the Cuntz
groupoid, $\eta_t(\Sab) = e^{it(\ell(\alpha)-\ell(\beta))}\Sab$.
In particular, for each
generating isometry $S_k$, $\eta_t(S_k)=e^{it}S_k$.  This determines
the automorphism $\eta_t$.  The dilation cocycle gives rise 
to the guage automorphisms which appeared in Theorem \ref{T:spec}.
\end{example}

The backward shift map on 
$X = \prod_1^{\infty}\{1,\dots,n\}$ is a useful tool
in the study of cocycles on $G_n$.
 Define $S \colon X \to X$  by
$S(x_1, x_2, x_3, \dots) = (x_2, x_3, \dots)$. 

Now let $(x,0, y) \in P_{UHF}$ and $k \geq 0$ be given.  
Then $(x, k, S^ky) \in G_n$.  Conversely,
if $(x, k, z) \in G_n$ with $k \geq 0$, there is a $y \in X$ such 
that $z = S^ky$ and $(x,0,y) \in P_{UHF}$.  
Indeed, we can take
$y = \alpha z$ where $\alpha$ is an arbitrary $k$-tuple.
If $k > 0$ and $(x,0,y) \in P_{UHF}$, 
then $(S^k x, -k, y) \in G_n$.  If 
$(z,-k,y) \in G_n$ with $k>0$, then there is $x \in X$ such that
$z = S^k x$ and $(x,0,y) \in G_n$.

\begin{theorem} \label{T:cocyc_fun}
If $d$ is a cocycle on $G_n$, define a continuous real valued
function $f_d$ on $X$ by
\[
f_d(x) = d(x,1,Sx), \quad \text{for all } x\in X.
\]

If $f \colon X \to \bbR$ is continuous, define a
cocycle $d_f$ on $G_n$ by
\begin{align}
d_f(x,k,S^ky) &= \sum_{j=0}^{k-1} f(S^jx) +
\sum_{j=k}^{\infty}[f(S^jx) - f(S^jy)], \label{pos} \\
d_f(S^k x,-k, y) &= -\sum_{j=0}^{k-1}f(S^jy) + 
\sum_{j=k}^{\infty}[f(S^jx) - f(S^jy)], \label{neg}
\end{align}
where $k \geq 0$ and $(x,0,y) \in P_{UHF}$.

These two correspondences are inverse to one another. Thus,
there is a one-to-one correspondence between continuous cocycles 
on $G_n$ and continuous real valued functions on $X$.
\end{theorem}

\begin{proof}
The continuity of $f_d$ follows from the continuity of $d$ and of the
map $x \mapsto (x,1,Sx)$.

Now suppose that $f$ is a continuous real valued map on $X$.  By the
remarks above, every element of $G_n$ has one of the two forms
$(x,k,S^k y)$ or $(S^k x, -k,y)$ where $k \geq 0$ and
$(x,0,y) \in P_{UHF}$.  In equation (\ref{pos}), the first $k$
coordinates of $y$ play no role in either side of the equation;
similarly, in equation (\ref{neg}) the first $k$ coordinates of $x$ play
no role.  When $k=0$, the two equations agree, since the first sum in
the right hand side is absent in either version.  Since
$(x,0,y) \in P_{UHF}$, $S^j x = S^j y$ for all large $j$; consequently
the sums in equations (\ref{pos}) and (\ref{neg}) are finite.

Observe that, with $k \geq 0$ and $(x,0,y) \in P_{UHF}$,
\begin{align*}
d_f (S^k y, -k, x) &= - \sum_{j=0}^{k-1} f(S^j x)
  +\sum_{j=k}^{\infty}[f(S^j y) - f(S^j x)] \\
&= - \left[ \sum_{j=0}^{k-1} f(S^j x) 
   +\sum_{j=k}^{\infty}[f(S^j x) - f(S^j y)] \right]\\
&= -d_f(x,k,S^ky).
\end{align*}
Now suppose that $k \geq 0$, $l \geq 0$,  $(x,0,y) \in P_{UHF}$ and
$(S^k y, 0, S^k z) \in P_{UHF}$.  It follows that
$(y,0,z) \in P_{UHF}$.  If we prove that
\[
d_f(x,k,S^k y) + d_f(S^k y, l, S^{k+l} z) =
d_f(x,k+l,S^{k+l} z).
\]
then the cocycle identity for $d_f$ on all of $G_n$
follows easily.

Now,
\begin{align*}
 &d_f(x,k,S^k y) + d_f(S^k y, l, S^{k+l} z) \\
&\qquad =  \sum_{j=0}^{k-1} f(S^jx) +
\sum_{j=k}^{\infty}[f(S^jx) - f(S^jy)]
+ \sum_{i=0}^{l-1} f(S^{i+k} y) 
+ \sum _{i=l}^{\infty}[f(S^{i+k} y) - f(S^{i+k} z)]\\
&\qquad = \sum_{j=0}^{k-1} f(S^j x) +
\sum_{j=k}^{k+l-1} f(S^j x) - \sum_{j=k}^{k+l-1} f(S^j y)
+\sum_{j=k+l}^{\infty} [f(S^j x) - f(S^j y)] \\
&\qquad \qquad +\sum_{j=k}^{k+l-1} f(S^{j} y) 
+\sum_{j=k+l}^{\infty}[f(S^j y) - f(S^j z)] \\
&\qquad = \sum_{j=0}^{k+l-1} f(S^j x) 
+ \sum_{j=l+l}^{\infty}[f(S^j x) - f(S^j z)] \\
&\qquad = d_f(x,k+l,S^{k+l}).
\end{align*}

It remains to show that these two maps are inverse to one another.

Let $f \colon X \to \bbR$ be continuous.  Then
\[
f_{d_f}(x) = d_f(x,1,Sx) =
f(x) + \sum_{j=1}^{\infty}[f(S^j x) - f(S^j x)] = f(x).
\]

It is a bit more complicated to show that the composition in the other
order is also the identity.

Let $d$ be a cocycle on $G_n$ and define $f\; (= f_d)$ by
$f(x) = d(x,1,Sx)$, for all $x \in X$.  If 
$(x,0,y) \in P_{UHF}$, then
\begin{align*}
f(x) + d(Sx,0,Sy) &= d(x,1,Sx) + d(Sx,0,Sy) \\
&= d(x,1,Sy) \\
&= d(x,0,y) + d(y,1,Sy) \\
&= d(x,0,y) + f(y)
\end{align*}
Therefore
\begin{equation}
f(x) - f(y) = d(x,0,y) - d(Sx,0,Sy). \label{diff}
\end{equation}
Apply this to $S^j x$ and $S^j y$ to obtain
\[
f(S^j x) - f(S^j y) = d(S^j x,0,S^j y) 
- d(S^{j+1} x,0,S^{j+1} y), \quad \text{for any } j \geq 0.
\]
Since $(x,0,y) \in P_{UHF}$, there is $N \in \bbN$ such that
$S^n x = S^n y$ (and hence $d(S^n x, 0, S^n y) = 0$) for all
$n \geq N$.  Addition yields
\begin{equation}
d(x,0,y) = \sum_{j=0}^{\infty}[f(S^j x) - f(S^j y)]. \label{diffsum}
\end{equation}
The sum is, of course, really a finite sum.

Next, we express $d(x,k,S^k y)$ in terms of $f$ (when $k > 0$).  Just
observe that
\begin{align*}
d(x,k,S^k y) &= d(x,1,Sx) + d(Sx,1,S^2 x) + \dots +
d(S^{k-1} x, 1, S^k x) + d(S^k x, 0, S^k y) \\
&= f(x) + f(Sx) + \dots + f(S^{k-1} x) + d(S^k x, 0, S^k y) \\
&= \sum_{j=0}^{k-1} f(S^j x) + 
\sum_{j=k}^{\infty}[f(S^j x) - f(S^j y)],
\end{align*}
where equation (\ref{diffsum}) is used in the last equality.

We have now shown that $d$ and $d_{f_d}$ agree on $(x,k,y)$ whenever
$k \geq 0$.  But $d$ and $d_f$ are both cocycles, so they also agree
whenever $k < 0$, i.e., on all of $G_n$. 
\end{proof}

The correspondence between continuous functions on $X$ and continuous
cocycles on $G_n$ is clearly a linear one.  It also preserves uniform
convergence on compacta.

\begin{proposition} \label{P:conv_uniform}
Let $f$ and $f_n$ ($n \in \bbN$) be continuous functions on $X$.
Then $f_n \to f$ uniformly on $X$ if, and only if,
$d_{f_n} \to d_f$ uniformly on compacta in $G_n$.
\end{proposition}

\begin{proof}
Assume that $f_n \to f$ uniformly on $X$.  Let 
$k \geq 0$ and let $\alpha$ and $\beta$ be such that
$\ell(\alpha) - \ell(\beta) =k$.  Let $(x,k,S^k y) \in \Uab$.
Then
\begin{align*}
&d_{f_n}(x,k,S^k y) - d_f(x,k,S^k y) = \\
&\qquad \sum_{j=0}^{k-1}[f_n(S^j x) - f(S^j x)] +
\sum_{j=k}^{\ell(\alpha)}\left([f_n(S^j x) - f(S^j x)]
- [f_n(S^j y) - f(S^j y)]\right)
\end{align*}
and this converges uniformly to 0 on $\Uab$.
A similar argument applies when
$\ell(\alpha) - \ell(\beta) < 0$.

For the converse, observe that
\[
f(x) - f_n(x) = d_f(s,1,Sx) - d_{f_n}(x,1,Sx)
\]
and that
\[
\{(x,1,Sx) \mid x \in X \} = \bigcup_{i=1}^n U_{(i),\emptyset}.
\]
Since $d_{f_n} \to d_f$ uniformly on this set,
$f_n \to f$ uniformly on $X$.
\end{proof}

The following remark will be useful in Theorems \ref{T:finrange} and
\ref{T:nococycles}

\begin{remark} \label{r:finval}
If $A$ is a compact and open subset of $X$, then $A$ is a finite
disjoint union of cylinder sets (sets of the form 
$E_{\alpha}=\{x\in X \mid x_1=\alpha_1, \dots, x_k=\alpha_k\}$).
Indeed, for each $a\in A$ there is a cylinder set $C_a$ such
that $a \in C_a \subseteq A$.  Since $A$ is compact, the open cover 
$\{C_a \mid a \in A\}$ has a finite subcover.  Thus $A$ is a finite
union of cylinder sets.  Since any cylinder set $E_{\alpha}$ can be
written as a disjoint union of cylinder sets of the form $E_{\beta}$
where the $\beta$'s all have a specified 
length greater than the length
of $\alpha$, $A$ can be written as a union of disjoint cylinder sets.

Now suppose that $f$ is a continuous function on $X$ which assumes
only finitely many values.  For each $t$
in the range of $f$, $f^{-1}(t)$ is
a compact open subset of $X$, and hence can be written as a finite
disjoint union of cylinder sets.  It follows that there is a positive
integer $N$ such that, for all $x$,
 $f(x)$ depends only on the first $N$ coordinates
of $x$.
\end{remark}

\begin{theorem} \label{T:finrange}
If $f$ is a continuous function on $X$ with finite range, then the
cocycle $d$ associated with $f$ is locally constant.  
\end{theorem}

\begin{proof}
By  Remark \ref{r:finval}, there is a positive
integer $N$ such that the value of $f$ at any point of $X$ depends
only on the first $N$ coordinates of the point.

Let $k \geq 0$ and $(x,0,y) \in P_{UHF}$.  We shall show that $d$ is
constant on a neighborhood of $(x,k,S^ky)$.  Let $P$ be an integer
greater than $k$ such that $x_i = y_i$ for all $i \geq P$.  Let
\begin{align*}
\alpha &= (x_1, \dots, x_{P+N}), \\
\beta &= (y_{k+1}, \dots, y_{P+N}).
\end{align*}
Then $\ell(\alpha) - \ell(\beta) = k$ and
$(x,k,S^k y) \in \Uab$.  For any $\gamma \in X$,
\[
d(\alpha \gamma, k, \beta \gamma) =
\sum_{j=0}^{k-1}f(S^j(\alpha \gamma)) +
\sum_{j=k}^P [f(S^j(\alpha \gamma)) - f(S^{j-k}(\beta \gamma))].
\]
With $j \leq P$, the first $N$ coordinates of $S^j(\alpha \gamma)$ are
all coordinates of $\alpha$ and the first $N$ coordinates of
$S^{j-k}(\beta \gamma)$ are all coordinates of $\beta$; thus the value
of $d$ is independent of $\gamma$.  This shows that
 $d$ is constant on
$\Uab$.  For points of the form $(S^k x, -k, y)$, choose a basic open
neighborhood $\Uab$ for $(y,k,S^k x)$ on which $d$ is constant: then 
$U_{\beta, \alpha}$ is a neighborhood of $(S^k x, -k, y)$
on which $d$ is constant.
Thus, $d$ is locally constant.
\end{proof}

\begin{corollary}
If $f$ is a continuous function on $X$ with finite range and $d$ is
the corresponding cocycle, then $d^{-1}(0)$ is open.  Consequently, 
$d^{-1}[0,\infty)$ is open and the support set of an analytic algebra.
\end{corollary}

\begin{proof}
Let $(x,k,z) \in G_n$ be such that $d(x,k,z)=0$.  Let $\Uab$ be a
neighborhood of $(x,k,z)$ on which $d$ is constant.  Then
$(x,k,z) \in \Uab \subseteq d^{-1}(0)$ and  $d^{-1}(0)$ is open.
\end{proof}

The dilation cocycle $d(x,k,y) = k$ is non-negative on the set
\begin{align*}
P_+ &= \{(\alpha\gamma, \ell(\alpha)-\ell(\beta), \beta\gamma 
\mid \gamma \in X, \ell(\alpha) \geq \ell(\beta) \}\\
&= \bigcup \{\Uab \mid \ell(\alpha) \geq \ell(\beta)\}.
\end{align*}
This is an open (and closed) subset of $G_n$ which is clearly closed
under the groupoid multiplication.  The analytic algebra associated
with the dilation cocycle is $A(P_+)$.  Since 
$P^{\vphantom{+}}_+ \cup P^{-1}_+ = G_n$
and $P^{\vphantom{+}}_+ \cap P^{-1}_+ = P_{UHF}$,
it follows that $A(P_+) + A(P_+)^*$ is dense in $O_n$ and
$A(P_+) \cap A(P_+)^* = A(P_{UHF})$, the canonical UHF subalgebra.

The diagonal algebra in $A(P_+)$ is not abelian, but
a modification gives  examples of  strongly maximal
triangular subalgebras of $O_n$.  Let
$Q$ be an open subset of $P_{UHF}$ which satisfies the following
\begin{enumerate}
\item $Q \circ Q \subseteq Q$
\item $Q \cap Q^{-1} = D_0$
\item $Q \cup Q^{-1} = P_{UHF}$
\end{enumerate}
Then $A(Q)$ is a strongly maximal triangular 
subalgebra of the canonical UHF
\cstar algebra, $A(P_{UHF})$.  Two specific examples are:
\begin{align*}
Q^{\text{ref}} &= \{(\alpha\gamma,0,\beta\gamma) \in P_{UHF} \mid \alpha
 \preceq \beta \} \\
Q^{\text{st}} &= \{(\alpha\gamma,0,\beta\gamma) \in P_{UHF} \mid \alpha
 \preceq_{\text{r}} \beta \} 
\end{align*}
  Here $\preceq$ is the lexicographic order and $\preceq_{\text{r}}$
  is the reverse lexicographic order.
 $Q^{\text{ref}}$ is the
spectrum of the refinement subalgebra of the canonical UHF algebra.
and $Q^{\text{st}}$ is the spectrum of the standard subalgebra.

Given $Q$ as above,
 let
\[
Q_{+} = \{(x,k,y) \mid \;\text{either } k=0 \text{ and }
(x,0,y) \in Q \text{ or } k >0 \}
 \]
It is easy to see that $Q_+ \circ Q_+ \subseteq Q_+$,
$Q_+ \cup Q_+^{-1} = G_n$, and 
$Q_+ \cap Q_+ = D_0$.  $A(Q_+)$ is a strongly maximal triangular
subalgebra of $\sO_n$.  We shall call $A(Q_+)$ the 
\emph{contractive extension} of $A(Q)$.

\begin{proposition} \label{semisimple_rad}
The  ``contractive'' algebras, $A(P_+)$ and $A(Q_+)$, are
semisimple.
\end{proposition}

\begin{proof}
By Theorem~\ref{T:spec}, both $A(P_+)$ and $A(Q_+)$ are invariant
under the guage automorphisms.  Therefore the radical is invariant
under guage automorphisms and so is generated by the Cuntz partial
isometries which it contains.  But if $\Sab$ is in the radical, then
so is  $S_\alpha = S_{\alpha}^{\phantom{*}} S_{\beta}^*
S_{\beta}^{\phantom{*}}$ (since $S_{\beta}$ is in $A(P_+)$ or
$A(Q_+)$, as appropriate).  Thus, if the radical contains a non-zero
operator, it contains one which is not quasi-nilpotent.  This is
impossible, so each algebra is semisimple.
 \end{proof}

\begin{example}
Let $f \colon X \to \bbR$ be given by $f(x)=1$ for all $x$.  Then it
is easy to see that
$d_f(x,k,z) = k$ for all $(x,k,z) \in G_n$.  Thus, the constant
function 1 corresponds to the dilation cocycle; the associated
analytic subalgebra is the contractive algebra $A(P_{+})$.

If $f(x)=c$ for all $x$ (with $c \neq 0$), then
$d_f(x,k,z) = ck$.  Again, the analytic 
subalgebra is $A(P_{+})$.
\end{example}

\begin{example}
Let $f \colon X \to \bbR$ by $f(x)=x_1$.  The range of $f$ is 
$\{1, \dots, n\}$.  When $k \geq 0$ and $(x,0,y)\in P_{UHF}$,
\[
d_f(x,k,S^ky) = x_1 + \dots + x_k + 
\sum_{i=k+1}^{\infty} (x_i - y_i)
\]
and
\[
d_f(S^k x, -k, y) = -y_1 - \dots - y_k +
\sum_{i=k+1}^{\infty} (x_i - y_i)
\]
For every $k$, there exist $x,z \in X$ such that $d_f(x,k,z)=0$;
i.e.~$d^{-1}(0)$ intersects every $k$-level set.  In particular,
$d^{-1}(0)$ properly contains $D_0$.
\end{example}

In UHF algebras, the $\bbZ$-analytic subalgebras are those analytic
algebras which arise from $\bbZ$-valued cocycles which vanish
precisely on the canonical diagonal.  Other important analytic
algebras, for example the refinement algebras, arise from bounded
cocycles which vanish precisely on the canonical diagonal.  Theorem
\ref{T:nococycles}
shows that the Cuntz algebras lack precise analogs of these
analytic subalgebras of UHF algebras.

\begin{theorem} \label{T:nococycles}
There is no cocycle $d$ on the Cuntz groupoid $G_n$ which is bounded
and which vanishes precisely on $D_0$  There is no cocycle $d$ which
is integer valued and which vanishes precisely on $D_0$.
\end{theorem}

\begin{proof}
First suppose that $d$ is a bounded cocycle.  Let $f$ be the
associated function on $X$ and let $\ol{1}=(1,1,\dots)$.  Then
$S\ol{1}=\ol{1}$.  If $f(\ol{1}) =0$, then 
$d^{-1}(0)$ properly contains $D_0$.  If
$f(\ol{1}) \neq 0$, then, for all $k>0$,
\[
d(\ol{1},k,S^k \ol{1}) = \sum_{j=0}^{k-1} f(S^j \ol{1}) =
\sum_{j=0}^{k-1} f(\ol{1}) = kf(\ol{1})
\]
and $d$ is unbounded.

Now suppose that $d$ is integer valued.  Then $f$ is also integer
valued.  Since $f$ is continuous and $X$ is compact, the range of $f$ 
is finite.  By  remark \ref{r:finval},
the value of $f$ at $x$
depends only on the first $N$ coordinates of $x$.

We will exhibit $x \neq y$ such that $(x,0,y) \in G_n$ and
$d(x,0,y)=0$.  Let $\alpha$ be the string consisting of
$2\cdots 2$ repeated $N-1$ times and let $\gamma$ be the infinite
string consisting entirely of 1's.  Set
\[
x=\alpha 21 \alpha \gamma \qquad \text{and} \qquad
y=\alpha 12 \alpha \gamma.
\]
It is routine to show that a string of length $N$ appears as the first
$N$ coordinates of one of the points
$x,Sx,\dots,S^{2N-1}x$ if, and only if, it appears as the first
$N$ coordinates of one of
$y,Sy,\dots,S^{2N-1}y$; furthermore, the frequency
of appearance is the same in both cases.
  Thus
\[
\sum_{j=0}^{2N-1} f(S^j x) = \sum_{j=0}^{2N-1} f(S^j y).
\]
Since $S^j x = S^j y$ for $j \geq 2N$, we have
\[
d(x,0,y) = \sum_{j=0}^{\infty}[f(S^j x) - f(S^j y)] = 0.
\]
Thus, $d^{-1}(0)$ properly contains $D_0$.
\end{proof}

As we pointed out earlier,
the standard TUHF algebra 
 is the strongly maximal triangular subalgebra of
$A(P_{UHF})$ whose spectrum is
$Q^{\text{st}}=\{(x,0,y) \in P_{UHF} \mid x \preceq_r y \}$.  
(The order is the
reverse lexicographic order.)  This is the prototypical 
$\bbZ$-analytic algebra.  The associated cocycle -- the standard
cocycle -- is given by the formula
\[
d(x,0,y) = \sum_{i=1}^{\infty} n^{i-1}(x_i-y_i).
\]

\begin{proposition} \label{P:stand_ext}
The standard cocycle $d$ on
$P_{UHF}$ has no extension to a cocycle on $G_n$.  
\end{proposition}

\begin{proof}
Suppose that $d$ has an extension to $G_n$, which we also denote by
$d$, and that $f$ is the associated function on $X$.  
Since 
\[
d(x,0,y) = \sum_{j=0}^{\infty}[f(S^j x) - f(S^j y)],
\]
any function which differs from $f$ by an additive constant will yield
a cocycle which agrees with the standard cocycle on $P_{UHF}$.
Therefore, we may, without loss of generality, assume that
$f(\ol{1})=0$, where $\ol{1}=(1,1,1,\dots)$.  Now let
$x=(2,2,\dots,2,\ol{1})$, where there are $N$ coordinates in $x$
which are 2.  By equation \ref{diff} in the proof of 
Theorem~\ref{T:cocyc_fun},
\[
f(x) = d(x,0,\ol{1})-d(Sx,0,S\ol{1})
= \sum_{j=1}^N n^{j-1} - \sum_{j=1}^{N-1} n^{j-1} = n^{N-1}
\]
This shows that $f$ is unbounded, an impossibility for a continuous
function on a compact set.  Thus the standard cocycle on $P_{UHF}$
admits no extension to $G_n$.
\end{proof}

A cocycle $d$ on $P_{UHF}$ is \emph{trivial}
if there is a continuous
function $b$ on $X$ such that
$d(x,0,y) = b(y) - b(x)$.  The refinement cocycle, given by
\[
d(x,0,y) = \sum_{i=1}^{\infty} \frac 1{n^i} (x_i - y_i),
\]
is one of the simplest and most important examples of a trivial
cocycle.  The function $b$ is given by
\[
b(x) = \sum_{i=1}^{\infty} \frac {x_i -1}{n^i};
\]
the range of $b$ is  $[0,1]$.  The
analytic subalgebra of $A(P_{UHF})$ associated with the refinement
cocycle is the refinement TUHF algebra; as we noted earlier,
the spectrum of this algebra is
$Q^{\text{ref}} = \{(x,0,y)\in P_{UHF} \mid x \preceq y \}$ 
(lexicographic order).  In
general, an analytic subalgebra of $A(P_{UHF})$ whose cocycle is
trivial is said to be trivially analytic.

 \begin{theorem} \label{T:triv_ext}
Let $A(Q)$ be a trivially analytic subalgebra of $A(P_{UHF})$
with cocycle $d$.  
Then $d$ has extensions to cocycles on $G_n$.  Furthermore,
there is  an extension so that the corresponding analytic
subalgebra of $\sO_n$ is $A(Q_+)$.
\end{theorem}

\begin{remark}
As noted above,
the refinement subalgebra of $A(P_{UHF})$ is a trivially analytic
subalgebra.  Thus $A(Q^{\text{ref}}_+)$ is an 
analytic subalgebra of $\sO_n$
which is an ``extension'' of the refinement TUHF algebra.
\end{remark}

\begin{proof}
Let $b$ be a continuous function on $X$, let $d$ be the cocycle
on $P_{UHF}$ given by $d(x,0,y)=b(y)-b(x)$ and let $A(Q)$ be the
analytic subalgebra of $A(P_{UHF})$ associated with the cocycle.
  
Let $X_0$ be the
equivalence class of the point $\ol{1} = (1,1,1,\dots)$.  
This is a dense set
in $X$. We wish to define a function $f$ on $X$ which will yield a
cocycle on $G_n$ which extends $d$.

Begin by setting $f(\ol 1)=0$. For $x \in X_0$, define
\[
f(x) = d(x,0,\ol 1)-d(Sx,0, \ol 1).
\]
Observe that, for any for any $x,y \in X_0$,
\begin{align*}
f(x)-f(y) &= d(x,0, \ol 1)-d(Sx, 0,\ol 1) 
- d(y,0,\ol 1) +d(Sy,0,\ol 1) \\
&= d(x,0, \ol 1)+ d(\ol 1, 0, y) 
-[d(Sx,0, \ol 1) + d(\ol 1, 0, Sy)]\\
&= d(x,0,y) - d(Sx,0,Sy)
\end{align*}

We claim that $f$ is uniformly continuous, and therefore has a
continuous extension to $X$.  Let $\epsilon >0$.
Let $\rho$ be the metric on $X$ defined by
\[
\rho(x,y) = \sum_{i=1}^{\infty}\frac {|y_i-x_i|}{n^i}.
\]
Since $b$ is uniformly continuous, there is $\delta >0$ such that
$|b(y) - b(x)| < \epsilon/2$ whenever $\rho(x,y) < \delta$.

Sets of the form $U_{\alpha} = \{\alpha\gamma \mid \gamma \in X\}$
form a neighborhood basis for the topology on $X$.  If the length of
$\alpha$ is $L$, then $\rho(x,y) \leq 1/n^L$ whenever
$x$ and $y$ are in $U_{\alpha}$.  Choose $L$ so that
$1/n^L < \delta$.  Let $U_{\alpha}$ be any basic open neighborhood
with $\ell({\alpha}) \geq L+1$.  Then, when
$x,y \in U_{\alpha} \cap X_0$, we have $Sx, Sy \in U_{S\alpha}$
and $\ell({S\alpha}) \geq L$.  Therefore,
\begin{align*}
|f(x)-f(y)| &= |d(x,0,y)-d(Sx,0,Sy)|
=|b(y)-b(x) -(b(Sy)-b(Sx))| \\
&\leq|b(y)-b(x)|+|b(Sy)-b(Sx)| < \epsilon.
\end{align*}

Denote the extension to $X$ by $f$ also.   We next show that the
relation 
\[
f(x) - f(y) = d(x,0,y) - d(Sx, 0, Sy)
\]
holds for any $(x,0,y) \in P_{UHF}$.  Let $k$ be such that 
$x_i = y_i$ for all $i \geq k$.  Let 
$\alpha = (x_1, \dots, x_k)$ and $\beta = (y_1, \dots, y_k)$.
Then $\Uab$ is a neighborhood of $(x,0,y)$ and a $G$-set in $G_n$.  In
particular, the projection maps onto the first and third coordinates
are each homeomorphisms when restricted to $\Uab$.  Let $x_{\nu}$ and
$y_{\nu}$ be sequences in $X_0$ such that
 $(x_{\nu}, 0 , y_{\nu}) \in \Uab$, $x_{\nu} \to x$ and 
$y_{\nu} \to y$.  It follows that
\begin{align*}
&d(x_{\nu},0,y_{\nu}) \to d(x,0,y), \\
&d(Sx_{\nu},0,Sy_{\nu}) \to d(Sx,0,Sy), \\
&f(x_{\nu}) \to f(x), \text{ and} \\
&f(y_{\nu}) \to f(y).
\end{align*}
Since
\[
f(x_{\nu}) - f(y_{\nu}) = 
d(x_{\nu},0,y_{\nu}) - d(Sx_{\nu}, 0, Sy_{\nu})
\]
for every $\nu$, we have
\[
f(x) - f(y) = d(x,0,y) - d(Sx, 0, Sy)
\]
as desired.

By an argument used in Theorem \ref{T:cocyc_fun}, 
\[
d(x,0,y) = \sum_{j=0}^{\infty}[f(S^j x) - f(S^j y)].
\]
It now follows immediately that the cocycle $d_f$
 on $G_n$ induced by the
function $f$ extends the cocycle $d$ on $P_{UHF}$.

Fix a value for $k$ and let 
$G_{n,k} = \{(x,j,y)\in G_n \mid j=k\}$.
We claim that $d_f$ is bounded on each $G_{n,k}$.
Indeed, when $k \geq 0$,
any point in $G_{n,k}$ can be written in the
form $(x,k,S^k y)$ for some $(x,0,y) \in P_{UHF}$ and
\begin{align*}
d_f(x,k,S^k y) &=  \sum_{j=0}^{k-1} f(S^j x) + 
\sum_{j=k}^{\infty}[f(S^j x) - f(S^j y)] \\
&=\sum_{j=0}^{k-1} f(S^j x) + d(S^k x,0,S^k y).
\end{align*}
Since $f$ is bounded, $k$ is fixed, and $d$ is bounded on
$P_{UHF}$, $d_f$ is bounded on $G_{n,k}$.
The boundedness of $d_f$ on $G_{n,k}$ when $k<0$ follows from the
cocycle property.

Finally, we show that there is an extension of $d$ such that the
corresponding analytic subalgebra is $A(Q_+)$.  Indeed, if $c$ is any
constant and $g = f + c$, then, since the correspondence between
functions on $X$ and cocycles is linear,
$d_g(x,k,y) = d_f(x,k,y) + kc$.  In particular,
$d_g(x,1,y) = d_f(x,1,y) + c$.  Since $d_f$ is bounded on
$G_{n,1}$, we may pick $c$ so that $d_g > 0$ on $G_{n,1}$.
For $k>0$, any element of $G_{n.k}$ is a product of $k$ elements of 
$G_{n,1}$; the cocycle property implies that $d_g > 0$ on $G_{n,k}$.
The cocycle property also guarantees that $d_g$ is negative on 
$G_{n,k}$ when $k<0$.  Thus, $d_g \geq 0$ precisely on $Q_+$.
 \end{proof}

\section{The Volterra Subalgebra of $\sO_n$} \label{S:Volterra}

The intersection of the Cuntz algebra  
with the Volterra nest algebra on $L^2[0,1]$ has been studied by
S.~C.~Power in~\cite{MR86d:47057}.  We refer to this intersection as
the \emph{Volterra subalgebra} of $O_n$.  
In this section, we show that
groupoid techniques allow us to obtain many of Power's results easily,
as well as some new information.
We first identify an open set $P_V$ for which $A(P_V)$ will turn out
to be the Volterra subalgebra of $O_n$.
An element $(\alpha\gamma, \ell(\alpha)-\ell(\beta), \beta\gamma)$ of
$G_n$ is in $P_V$ if any of the following five conditions
on $\alpha$ and $\beta$ hold:
\begin{gather}
\ell(\alpha)=\ell(\beta) \text{ and } \alpha \preceq \beta,
 \text{ or} \label{Veq1}  \\
\ell(\alpha) < \ell(\beta) \text{ and }
\alpha \prec
(\beta_1, \dots, \beta_{\ell(\alpha)}), \text{ or} \label{Veq2}\\
\ell(\alpha) < \ell(\beta) \text{ and }
\alpha =
(\beta_1, \dots, \beta_{\ell(\alpha)})
 \text{ and }
\beta_{\ell(\alpha)+1}= \dots = \beta_{\ell(\beta)}=n,
  \text{ or}  \label{Veq3}\\
\ell(\alpha) > \ell(\beta) \text{ and }
(\alpha_1, \dots, \alpha_{\ell(\beta)}) \prec
\beta, \text{ or} \label{Veq4}\\
\ell(\alpha) > \ell(\beta) \text{ and }
(\alpha_1, \dots, \alpha_{\ell(\beta)}) =
\beta \text{ and }
\alpha_{\ell(\beta)+1}= \dots = \alpha_{\ell(\alpha)}=1. \label{Veq5}
\end{gather}
The order on finite strings used in these conditions is the
lexicographic order.

There is an alternate description of $P_V$.  Although the
lexicographic order is most commonly used for finite sequences or for
sequences in $X$ which are  equivalent in the sense that `tails are
equal', the lexicographic order is also useful as a total order on all
of $X$.  With this in mind, it is not hard to see that
$P_V$ is the union of the following four sets:

\begin{align*}
R &= \{(x,k,y) \mid x \prec y\} \\
D_0 &=\{(x,0,x) \mid x \in X \} \qquad \text{(diagonal)}\\
S_e &=\{(x,k,x) \mid k < 0 \text{ and } x_j = n \text{ for all large }
j \} \qquad \text{(expansive)}\\
S_c &= \{(x,k,x) \mid k > 0 \text{ and } x_j = 1 \text{ for all large
  } j \} \qquad \text{(contractive)}
\end{align*}

Indeed, if $x \prec y$ and $(x,k,y) \in G_n$, let $j$ be the first
index for which $x_j \neq y_j$.  Choose initial segments $\alpha$ and
$\beta$ of $x$ and $y$ such that $j < \min(\ell(\alpha),\ell(\beta))$
  and $k = \ell(\alpha) - \ell(\beta)$.  Then $(x,k,y)$ has the form
$(\alpha\gamma,\ell(\alpha)-\ell(\beta),\beta\gamma)$ for suitable
$\gamma$ and also satisfies one of conditions (\ref{Veq1}),
(\ref{Veq2}), or (\ref{Veq4}).  Elements of the form $(x,0,x)$ clearly
satisfy condition (\ref{Veq1}).
If $k<0$ and $x_j = n$ for all large $j$, then $(x,k,x)$ satisfies
condition (\ref{Veq3}); if $k>0$ and $x_j = 1$ for all large $j$, then
$(x,k,x)$ satisfies condition (\ref{Veq5}).

In the other direction, if $(x,k,y)$ satisfies condition (\ref{Veq2})
or condition (\ref{Veq4}), then $x \prec y$.
Suppose that 
$(x,k,y) = (\alpha\gamma, \ell(\alpha) - \ell(\beta), \beta\gamma)$
 satisfies (\ref{Veq3}).  If  any of
$\gamma_1, \dots, \gamma_{\ell(\beta)-\ell(\alpha)}$ are less than $n$,
then $x \prec y$.  If all of
$\gamma_1, \dots, \gamma_{\ell(\beta)-\ell(\alpha)}$ equal $n$ and any
 $\gamma_j < n$ with $j$ between $\ell(\alpha)-\ell(\beta)+1$ and
$2(\ell(\beta)-\ell(\alpha))$ then again $x \prec y$.  Continuing in
this fashion, we see that either $x \prec y$ or $x=y$ and $x_j=n$ for
all large $j$.  Similarly, if $(x,k,y)$ satisfies condition
(\ref{Veq5}), then either $x \prec y$ or $x_j =1$ for all large $j$.
Finally, if $(x,k,y)$ satisfies (\ref{Veq1}), then $k=0$ and $(x,k,y)$
is in $D_0 \cup R$.

Next we show that $P_V$ is an open subset of $G_n$.  Since $D_0$ and
$\{(x,k,y) \in G_n \mid x \prec y \}$ are open sets, we need only show
that all points $(x,k,x)$
 with $k<0$ and a tail consisting of $n$'s or with
$k>0$ and a tail consisting of $1$'s lie in open sets contained in
$P_V$.  Suppose that $k>0$ and $x=\beta\delta$, where
$\beta = (\beta_1, \dots, \beta_s)$ and $\delta = (1,1,1, \dots )$.
Let $\alpha = (\beta_1, \dots, \beta_s, 1, \dots, 1)$, where there are
$k$ 1's following $\beta_s$.  Then
$U_{\alpha,\beta} = \{(\alpha\gamma, k, \beta\gamma) \mid \gamma \in
X\}$ is an open subset of $G_n$.  If $\gamma=(1,1,\dots)$, then
$(\alpha\gamma, k, \beta\gamma)=(x,k,x)$,
so $(x,k,x) \in U_{\alpha,\beta}$.
If $\gamma \neq (1,1,\dots)$, then it is easy to see that
$\alpha\gamma \prec \beta\gamma$ and hence
$(\alpha\gamma, k, \beta\gamma) \in P_V$.  Thus
$(x,k,x) \in  U_{\alpha,\beta} \subseteq P_V$.  
A similar argument takes
care of points with  $k<0$ and a tail consisting of $n$'s; thus
$P_V$ is open.

While $P_V$ is an open subset of $G_n$, it is not closed.  To see
this, let $\alpha$ be any finite string and let $\delta$ be a string
of length $k$ in which some $\delta_i \neq 1$.  Since $\delta$ is not
first in the lexicographic order on
$\{1, \dots, n\}^k$, there is a string $\eta$ of length $k$ such that
$\eta \prec \delta$.  For each $p \in \bbN$, let
\begin{align*}
&y^p = \alpha\delta \dots \delta\eta \eta \dots
&& (p \text{ copies of }\delta), \\
&z^p = \alpha\delta\delta \dots \delta\eta\eta \dots
&& (p+1 \text{ copies of }\delta), \\
&x = \alpha \delta \delta \dots && (\text{infinitely many
copies of } \delta).
\end{align*}
Since $y^p \prec z^p$, $(y^p, k, z^p) \in P_V$.  Also
$(y^p, k, z^p) \longrightarrow (x,k,x)$ in the topology on $G_n$.
But $(x,k,x) \notin P_V$, so $P_V$ is not closed.

When is a partial isometry $\Sab$ in $A(P_V)$?
For the case when $\ell(\alpha)=\ell(\beta)$, i.e., when
$\Sab$ is in the canonical UHF subalgebra,
$\Sab$ is  in $A(P_V)$ if, and only if,
$\Sab$ is also in the Volterra nest
subalgebra on $L^2[0,1]$.
These partial isometries correspond to
condition \ref{Veq1}.
Expansive partial isometries in $A(P_V)$
correspond to conditions
\ref{Veq2} and \ref{Veq3}.  Condition \ref{Veq2} gives the expansive
partial isometries for which the range interval lies to the left of
the domain interval.  Condition \ref{Veq3} gives those expansive
partial isometries for which the domain interval is a subinterval
of the range interval
located at the right end of the range interval.  Similarly, conditions
\ref{Veq4} and \ref{Veq5} correspond to the contractive partial
isometries in $A(P_V)$.
  Condition \ref{Veq4} gives the contractive partial
isometries for which the range interval  lies to the left of the
domain interval and condition \ref{Veq5} describes
 the case when the range
interval is a subinterval of the domain interval which is located at
the left end of the domain interval.

From these considerations, we see that
the partial isometries $\Sab$
which are in $A(P_V)$ are precisely the 
Cuntz partial isometries 
 which lie in the Volterra nest algebra.
Since $A(P_V)$ is generated by its Cuntz partial isometries, it
follows that $A(P_V)$ is contained in 
$O_n \cap \Alg \sV$, where $\sV$ is the Volterra nest.

  While most of
the projections in $\sV$ are not in $O_n$, those that are are
strongly dense in $\sV$.  These projections can be described in
several ways.  The most useful is as follows:  let
$T = \{x \in X \mid x_j = n \text { for all large } j\}$.
For $x \in T$, let $p_x$
be the characteristic function of
$\{(y,0,y) \mid y \preceq x\}$.
This subset of $P_V$ is compact and, since $x \in T$,
also open; thus
$p_x$ is a continuous function with compact support
and can therefore be viewed as an element of $A(P_V)$.

As an operator on $L^2[0,1]$ (using the representation of the Cuntz
algebra described above), $p_x$ is the projection onto the subspace
$L^2[0,x]$, where
\[
x = \sum_{j=1}^{\infty} \frac {x_j-1}{n^j}.
\]
And in terms of the generators of the Cuntz algebra, each $p_x$ has
the form
 $ \sum_{\alpha} S^{\vphantom{*}}_{\alpha}S^*_{\alpha}$, where
all $\alpha$ in the summation have the same length, say $k$, and run
 through an initial segment in the set of $k$-tuples with the
 lexicographic order.

\begin{theorem} \label{T:Volterra}
$A(P_V) = \sO_n \cap \Alg \sV$.
\end{theorem}

\begin{proof}
  If $f \in O_n$ and 
$x \in T  = \{x \in X \mid x_j = n \text { for all large } j\}$, 
then for any $(s,k,t) \in G_n$,
\begin{align*}
(f \ast p_x)(s,k,t) &=f(s,k,t)p_x(t,0,t) \\
 &=
\begin{cases}
f(s,k,t), \quad &\text{if } t \preceq x, \\
0, \quad &\text{if } x \prec t
\end{cases}
\end{align*}
and
\begin{align*}
(p_x \ast f \ast p_x) (s,k,t) &= p_x(s,0,s)f(s,k,t)p_x(t,0,t) \\
&=
\begin{cases}
f(s,k,t), \quad &\text{if both } s \preceq x \text{ and }
t \preceq x \\
0, \quad &\text{otherwise}.
\end{cases}
\end{align*}

Now suppose that $f \in O_n \cap \Alg \sV$.  
To prove that $f \in A(P_V)$
we must show that $f$ vanishes on the complement of $P_V$ in $G_n$.
First, suppose that $t \prec s$ and $(s,k,t) \in G_n$.
If $t$ is an immediate predecessor of $s$, choose $x=t$; otherwise
choose $x \in T$ such that $t \prec x \prec s$.  In either case,
$x \prec s$.   Since $f \in \Alg \sV$ and 
$p_x \in \sV$, we have
$f \ast p_x = p_x \ast f \ast p_x$. 
 It follows from the equations above that
$f(s,k,t) = 0$.  This leaves the case in which $s=t$.  There are two
possibilities: $k>0$ and infinitely many $t_i >1$, for which the
details are below, and $k<0$ and infinitely many $t_i<n$ (which can be
done in a similar fashion).

Suppose that $(t,k,t) \in G_n$ with $k>0$ and infinitely many
$t_i>1$. Then there are finite strings $\tau$ and $\gamma$ such that
one of the components of $\gamma$ is greater than 1,
$\ell(\gamma)=k$, and
$t = \tau \ol{\gamma}$.  (Here, $\ol{\gamma}$ denotes the
infinite concatenation $\gamma \gamma \gamma \dots$.)
If $f(t,k,t) \neq 0$, then there is a neighborhood $U$ of $(t,k,t)$
such that $f$ is non-zero on $U$.  For $n$ a positive integer, let
$\gamma^n$ denote $\gamma \dots \gamma$, the $n$-fold concatenation of
$\gamma$ with itself.  Then there exists $n$ such that if
$\alpha = \tau \gamma^{n+1}$ and $\beta = \tau \gamma^n$, then
$U_{\alpha,\beta} \subseteq U$.  Let $\delta$ be any sequence of
length $k$ which strictly precedes $\gamma$ in the lexicographic
order.  (The existence of such a sequence is guaranteed by the fact
that one of the components of $\gamma$ is greater than 1.)  Let
$x = \alpha\ol{\delta}$ and $y = \beta \ol{\delta}$.  Then
$y \prec x$ and $(x,k,y) \in U_{\alpha,\beta} \subseteq U$. By
the preceding paragraph,  $f(x,k,y) = 0$; but this contradicts the
fact that $f$ does not vanish on $U$.
Thus, if $k>0$ and infinitely many $t_i >1$, then $f(t,k,t)=0$.
Similarly, if $k<0$ and infinitely many $t_i < n$, then
$f(t,k,t)=0$.  Thus, $f$ vanishes on the complement of $P_V$ and we
have shown that $O_n \cap \Alg \sV \subseteq A(P_V)$. We
have already verified the reverse containment, so equality is 
proven.
\end{proof}

\begin{proposition}  \label{volterra_rad}
The spectrum of the radical in $A(P_V)$ is
$\{(x,k,y) \in G_n \mid x \prec y\}$.
\end{proposition}

\begin{proof}
Let $R=\{(x,k,y) \in G_n \mid x \prec y\}$ and let $R_0$ denote the
spectrum of $\rad A(P_V)$, the Jacobson radical of $A(P_V)$. 
 We first show that
$R \subseteq R_0$ while the three sets
\begin{align*}
&D_0 = \{(x,0,x) \mid x \in X \}, \\
&S_e =\{(x,k,x) \mid k<0, x \in X
\text{ and $x_j = n$, for all large $j$} \}, \\
&S_c =\{(x,k,x) \mid k>0, x \in X
\text{ and $x_j = 1$, for all large $j$} \}
\end{align*}
are all disjoint from $R_0$; in other words,
$R=R_0$.

Let $(x,k,y)\in R$.  Then there exist finite strings
$\alpha$ and $\beta$ and $\gamma \in X$ such that
$\ell(\alpha) - \ell(\beta) =k$, $\alpha \prec \beta$,
$x = \alpha \gamma$ and $y = \beta \gamma$.
The set $U_{\alpha,\beta}$, the support set for
$\Sab$, is a subset of $R$.
Since the ``range interval'' for the partial isometry
$\Sab$ lies to the left of the
``domain interval'', there is a projection $p \in \sV$
such that $\Sab = p \Sab p^{\perp}$.

It follows easily that, for any $T \in A(P_V)$, both
$\Sab T$ and $T\Sab$ are nilpotent
(of order 2).  Alternatively, since $\Sab$ lies in
the radical of the Volterra nest algebra, the ideal generated by
$\Sab$ in the Volterra nest algebra consists entirely
of quasi-nilpotents; a fortiori, the ideal generated by
$\Sab$ in $A(P_V)$ consists of quasi-nilpotents.
Either way, $\Sab$ lies in the radical of
$A(P_V)$.  Thus $U_{\alpha,\beta} \subseteq R_0$; in particular,
$(x,k,y) \in R_0$.  Thus $R \subseteq R_0$.

We next show that $D_0 \cap R_0 = \emptyset$.  Suppose that
$(x,0,x) \in D_0 \cap R_0$.  Then there is a finite string
$\alpha$ such that $U_{\alpha,\alpha} \subseteq R_0$.  But this
implies that $S^{\phantom{*}}_{\alpha}S^*_{\alpha}$, a 
non-zero projection,
lies in $\rad A(P_V)$.  This is a contradiction, so
$D_0 \cap R_0 = \emptyset$.

The final two pieces, $S_e \cap R_0 = \emptyset$ and
 $S_c \cap R_0 = \emptyset$, are proven in 
essentially the same way, so
 we provide only one of the proofs.
Suppose $(x,k,x) \in S_e \cap R_0$. 
 Since $S_e$ and $R_0$ are open, there
is a neighborhood of $(x,k,x)$ which is a subset of $R_0$.  It follows
that there are finite strings $\alpha$ and $\beta$ corresponding the
intervals $[a,c]$ and $[b,c]$ such that $a<b$ and
$\Sab \in \rad A(P_V)$.  But if we compress
$\Sab$ to $L^2[a,c]$, we get a co-isometry,
all powers of which have norm 1.  So all powers of
$\Sab$ have norm 1, contradicting
the assertion that $\Sab \in \rad A(P_V)$.  This shows that
$S_e \cap R_0 = \emptyset$.

In the proof that $R \subseteq R_0$, we showed that if
$(x,k,y) \in R$, then there is a neighborhood
$\Uab$ of $(x,k,y)$ such that $\Uab \subseteq R$ and
$\Sab \in \rad A(P_V)$.  We next show that this implies that
$A(R) \subseteq \rad A(P_V)$.  Since $A(R)$ is generated by the Cuntz
partial isometries which it contains, it is sufficient to show that each 
Cuntz partial isometry in $A(R)$ lies in the radical.

If $\Scd$ is a Cuntz partial isometry in $A(R)$, then for each
$(x,k,y) \in \Ucd$ there is a neighborhood $\Uab$ such that
$(x,k,y) \in \Uab \subseteq \Ucd$ and $\Sab \in \rad A(P_V)$.
This gives an open cover of $\Ucd$, from which we can select a finite
subcover.  It is routine to arrange that the sets in the subcover are
pairwise disjoint (with corresponding Cuntz partial isometries still
in the radical).  
  So $\Ucd$ can be written as a finite disjoint
union of sets of the form $\Uab$ for which $\Sab \in \rad A(P_V)$.
But then $\Scd$ is the sum of these $\Sab$, whence
$\Scd \in \rad A(P_V)$.  Thus $A(R) \subseteq \rad A(P_V)$.

We now know that $A(R) \subseteq \rad A(P_V) \subseteq A(R_0)$ and
that $R=R_0$.  But then $A(R) = A(R_0)$; in particular,
$\rad A(P_V) = A(R)$.
\end{proof}

\begin{remark}
If we let $B$ denote the norm closure 
of $A(D_0)+\rad A(P_V)$ in $A(P_V)$, then
$B$ is a proper subset of $A(P_V)$.  This follows from the fact that
the spectrum of $B$ is $D_0 \cup R$.  While far from being
semi-simple, the algebra $A(P_V)$ falls short of having a radical plus
diagonal decomposition.
\end{remark}

As mentioned earlier,
the  Volterra  subalgebra $A(P_V)$ has been studied
extensively by Power in~\cite{MR86d:47057}.  Making use of the
spectrum, we can obtain several of Power's results with new proofs
which provide a different intuitive insight.  For example, Power
points out that $A(P_V)$ is \emph{non-Dirichlet} in the sense that
$A(P_V) + A(P_V)^*$ is not dense in $O_n$; in other language, $A(P_V)$
is triangular but not strongly maximal triangular.  This is evident
from the fact that $P_V \cup P_V^{-1} \neq G_n$.

We prove the next theorem, which appears in~\cite{MR86d:47057}, using
spectral techniques.  In this theorem, $\com A(P_V)$ 
denotes the closed
ideal generated by all commutators in $A(P_V)$.

 \begin{theorem} \label{com_ideal}
The radical of $A(P_V)$ is equal to the commutator ideal of
$A(P_V)$.
\end{theorem}

\begin{proof}
We begin by showing that $\com A(P_V) \subseteq \rad A(P_V)$.  It
suffices to show that $[f,g] \in \rad A(P_V)$ for all $f$ and $g$ in
$A(P_V)$.  To do this, it is enough
 to prove that $[f,g]$ vanishes on
$D_0 \cup S_e \cup S_c$.  For any $(x,k,x) \in D_0 \cup S_e \cup S_c$,
the product $f \ast g$ is given by
\[
f \ast g (x,k, x) = \sum f(x,i,u) g(u, k-i, x)
\]
The sum is taken over all $i \in \bbZ$ and
 $u$ for which $(x,i,u)$ and $(u,k-i,x)$ lie in $P_V$.  But this
requires $x \preceq u$ and $u \preceq x$; thus $u=x$.  Furthermore,
$i$ and $k-i$ cannot have opposite signs; otherwise $x$ would have to
possess a tail consisting only of $n$'s and a tail consisting only of
$1$'s.  If $N_k$ is the set $\{0, \dots ,k\}$ when $k \geq 0$ and the
set $\{k, \dots, 0\}$ when $k \leq 0$, then the formula reduces to
\[
f \ast g (x,k, x) = \sum_{i \in N_k} f(x,i,x) g(x, k-i, x).
\]
The change of index $j=k-i$ now yields:
\begin{align*}
f \ast g (x,k, x) &= \sum_{i \in N_k} f(x,i,x) g(x, k-i, x)\\
 &= \sum_{j \in N_k} g(x,j,x) f(x, k-j,x) \\
 &= g \ast f (x,k,x).
\end{align*}
This shows that $f \ast g - g \ast f$ vanishes on
$D_0 \cup S_e \cup S_c$.  Thus $[f,g] \in A(R) = \rad A(P_V)$.

It remains to show that $\rad A(P_V) = A(R) \subseteq \com A(P_V)$.
By means of the same compactness argument used in the proof of
Proposition \ref{volterra_rad}, it suffices to show that for each
point $(x,k,y) \in R$, there is a neighborhood $\Uab$ of
$(x,k,y)$ contained in $R$ such that
$\Sab \in \com A(P_V)$.

Let $(x,k,y) \in R$.  Then $x \prec y$.  As in Proposition
\ref{volterra_rad}, there are finite strings $\alpha$ and $\beta$ and
an infinite string $\gamma$ such that $\alpha \prec \beta$,
$x = \alpha \gamma$, $y = \beta \gamma$, and
$k = \ell(\alpha) - \ell(\beta)$.  Furthermore, the partial isometry
$\Sab \in \rad A(P_V)$.  We know that there is a projection $p$ in the
Volterra nest such that
$\Sab = p \Sab p^{\perp}$.  Since the range interval for $\Sab$ lies
to the left of the domain interval, we can even arrange that $p$
corresponds to an initial segment in $[0,1]$ with $n$-adic right
endpoint; in particular, we may assume that $p \in A(P_V)$.  We then
have $[p, \Sab] = p\Sab - \Sab p = \Sab$.  Thus
$\Sab \in \com A(P_V)$.  
\end{proof}

\begin{remark}
The proof of the preceding theorem in~\cite{MR86d:47057} makes use of
 an intermediate result: the closed commutator ideal is equal to
the closed linear span of
$\{a \in O_n \mid a = p_s a (1-p_t) \text{ for some } s <
t\}$, which we denote by $J$.
Since the generators of $J$ are obviously in 
$\rad A(P_V) = \com A(P_V)$, we can establish this result in our
framework by showing that
$\rad A(P_V) \subseteq J$.  But if $(x,k,y) \in R$, then
(as above) there are finite strings $\alpha \prec \beta$ such
that $(x,k,y) \in \Uab \subseteq R$ and the range interval associated
 with $\Sab$ lies to the left of the domain interval associated with
$\Sab$.  Hence we can find $s<t$ such that
$\Sab = p_s \Sab (1-p_t)$.  Any Cuntz partial isometry in
$\rad A(P_V) = A(R)$ can be written as a sum of Cuntz partial
 isometries of this type, so $A(R) \subseteq J$.
\end{remark}

\begin{remark}
Power also proves that the (commutative) algebra
$A(P_V) / \rad A(P_V)$ is isomorphic to a function algebra.
The function algebra is defined on a subset $Y$ of
 $X \times \ol{\bbD}$;
this subset consists of all points $(x,0)$ for which $x$ is not
$n$-adic together with all points $(x,z)$ where $x$ is $n$-adic and
$z \in \ol{\bbD}$, the closed unit disk in $\bbC$.
We will not reprove this result, but it is worthwhile to state
explicitly the homomorphism from $A(P_V)$ (viewed as functions
supported on $P_V)$) to this function algebra which induces Power's
isomorphism on $A(P_V) / \rad A(P_V)$.

Let $S = D_0 \cup S_e \cup S_c = P_V \setminus R$ and let $\Phi$ be
the homomorphism on $A(P_V)$ which induces
Power's isomorphism.
For $f \in A(P_V)$,
$\Phi(f)$  depends only on $f|_S$.  If $x$ is not $n$-adic, then
$\Phi(f)(x,0) = f(x,0,x)$.  If $x$ is $n$-adic with a tail consisting
of $1$'s, then
\[
\Phi(f)(x,z) = \sum_{k=0}^{\infty} f(x,k,x) z^k.
\]
And, if $x$ is $n$-adic with a tail consisting of $n$'s, then
\[
\Phi(f)(x,z) = \sum_{k=-\infty}^0 f(x,k,x) z^{-k}
\]
When $f \in \rad A(P_V)$, $\Phi(f)=0$; therefore $\Phi$ induces a map
on the quotient of $A(P_V)$ by its radical.  This map is an isometric
isomorphism of operator algebras.
In particular, as in~\cite{MR86d:47057}, if $\alpha$ and $\beta$ are
strings corresponding to intervals which share a left endpoint or
share a right endpoint, and if $\Sab \in A(P_V)$, then
$\Phi(\Sab)$ is the monomial $z^k$, where $k$ is the absolute value of
the index of dilation for $\Sab$.
\end{remark}

We have seen how spectral techniques give new 
proofs for several of Power's results on the Volterra subalgebra of 
$O_n$.  The next proposition proposition contains new information
about the Volterra subalgebra.

\begin{proposition} \label{P:maximal}
$A(P_V)$ is a maximal triangular
subalgebra of $O_n$.
\end{proposition}

\begin{proof}

Suppose that $\sT$ is a triangular subalgebra of $O_n$ which
contains $A(P_V)$.  To prove that $A(P_V)$ is maximal triangular, we
need to show that $\sT \subseteq A(P_V)$.  To do this, it is
sufficient to show that if $f \in \sT$ and if $x$ is an $n$-adic
element of $[0,1]$, then $p_x^{\perp} f p_x = 0$.  It then follows
that $f$ leaves invariant every projection in the Volterra nest and so
lies in $O_n \cap \Alg \sV = A(P_V)$.

Observe that $p_x f p_x^{\perp}$ leaves invariant every projection in
the Volterra nest and that, when $x$ is $n$-adic, $p_x \in O_n$.  So
$p_x f p_x^{\perp} \in A(P_V) \subseteq \sT$.  Therefore,
$g = p_x^{\perp} f p_x + p_x f p_x^{\perp}$ is a self-adjoint element
of $\sT$.  Consequently, $g$ lies in the diagonal of $\sT$.  But the
diagonal is $A(D_0)$, so $g$ commutes with $p_x$.  Therefore
$0 = p_x^{\perp} g p_x = p_x^{\perp} f p_x$ and the maximal
triangularity follows.
\end{proof}

\begin{remark}
There is no
 analytic subalgebra $A(P)$ defined by a cocycle $d$ with the property
that 
$P_V \setminus D_0 \subset \{ (x, k, y) \in G_n
\mid d(x, k, y) > 0 \}$. 
Call this latter set $P^+$ and similarly define
$P^-$ and $P^0$.

Now $A(P^+ \cup D_0)$ is a triangular algebra 
containing $A(P_V)$;  by Proposition~\ref{P:maximal} the 
two are equal.  Write

\begin{align*}
G_n &= [(P^+ \cup D_0) \cup (P^- \cup D_0)] \cup [P^0 \setminus D_0] \\
    &= [P_V \cup P_V^{-1}] \cup [P^0 \setminus D_0].
\end{align*}

As both $P^0$ and $D_0$ are clopen, so is their difference.  
Since the two terms in square brackets are disjoint, we conclude
that $P_V \cup P_V^{-1}$ is closed.  But the same argument 
that was used to show that $P_V$ is not closed also shows that
$P_V \cup P_V^{-1}$ is not closed.  Thus we have a contradiction.
This shows that $A(P_V)$ is not contained in an analytic subalgebra of
$O_n$ with diagonal $A(D_0)$.
\end{remark}

\providecommand{\bysame}{\leavevmode\hbox to3em{\hrulefill}\thinspace}
\providecommand{\MR}{\relax\ifhmode\unskip\space\fi MR }
\providecommand{\MRhref}[2]{%
  \href{http://www.ams.org/mathscinet-getitem?mr=#1}{#2}
}
\providecommand{\href}[2]{#2}


\begin{thebibliography}{10}

\bibitem{MR57:7189}
Joachim Cuntz, \emph{Simple ${C}\sp*$-algebras generated by isometries}, Comm.
  Math. Phys. \textbf{57} (1977), no.~2, 173--185. \MR{57 \#7189}

\bibitem{MR82f:46073a}
Joachim Cuntz and Wolfgang Krieger, \emph{A class of {$C\sp{\ast} $}-algebras
  and topological {M}arkov chains}, Invent. Math. \textbf{56} (1980), no.~3,
  251--268. \MR{82f:46073a}

\bibitem{MR2002a:47107}
Kenneth~R. Davidson, Elias Katsoulis, and David~R. Pitts, \emph{The structure
  of free semigroup algebras}, J. Reine Angew. Math. \textbf{533} (2001),
  99--125. \MR{2002a:47107}

\bibitem{MR2000k:47005}
Kenneth~R. Davidson and David~R. Pitts, \emph{Invariant subspaces and
  hyper-reflexivity for free semigroup algebras}, Proc. London Math. Soc. (3)
  \textbf{78} (1999), no.~2, 401--430. \MR{2000k:47005}

\bibitem{MR98g:46083}
Alex Kumjian, David Pask, Iain Raeburn, and Jean Renault, \emph{Graphs,
  groupoids, and {C}untz-{K}rieger algebras}, J. Funct. Anal. \textbf{144}
  (1997), no.~2, 505--541. \MR{98g:46083}

\bibitem{MR90m:46098}
Paul~S. Muhly and Baruch Solel, \emph{Subalgebras of groupoid ${C}\sp
  *$-algebras}, J. Reine Angew. Math. \textbf{402} (1989), 41--75.
  \MR{90m:46098}

\bibitem{MR2001a:22003}
Alan L.~T. Paterson, \emph{Groupoids, inverse semigroups, and their operator
  algebras}, Birkh\"auser Boston Inc., Boston, MA, 1999. \MR{2001a:22003}

\bibitem{MR92k:47073}
Gelu Popescu, \emph{von {N}eumann inequality for {$(B({\mathcal H})\sp n)\sb
  1$}}, Math. Scand. \textbf{68} (1991), no.~2, 292--304. \MR{92k:47073}

\bibitem{MR86d:47057}
S.~C. Power, \emph{On ideals of nest subalgebras of ${C}\sp \ast$-algebras},
  Proc. London Math. Soc. (3) \textbf{50} (1985), no.~2, 314--332.
  \MR{86d:47057}

\bibitem{MR94g:46001}
Stephen~C. Power, \emph{Limit algebras: an introduction to subalgebras of
  {$C\sp *$}-algebras}, Pitman Research Notes in Mathematics Series, vol. 278,
  Longman Scientific \& Technical, Harlow, 1992. \MR{94g:46001}

\bibitem{MR82h:46075}
Jean Renault, \emph{A groupoid approach to ${C}\sp{\ast} $-algebras}, Springer,
  Berlin, 1980. \MR{82h:46075}
\end{thebibliography}
\end{document}